\documentclass{article}

\usepackage{longtable}
\usepackage[english]{babel} 
\usepackage{amsmath,amsfonts,amsthm,bm} 
\usepackage{graphicx}
\usepackage{gensymb}
\usepackage{dcolumn}
\usepackage{bm}
\usepackage{color}
\usepackage[autostyle]{csquotes}
\usepackage{tikz}
\usepackage{setspace}
\usepackage{multirow}
\usepackage{authblk}

\usepackage[latin1]{inputenc}
\usepackage{hyperref}

\newtheorem{conjecture}{Conjecture}

\theoremstyle{definition}

\usepackage{appendix}
\usepackage{enumerate}
\usepackage{listings}
\usepackage{tikz}
\usepackage{gensymb}
\usepackage{dirtytalk}

\usepackage{geometry}
\begin{document}

\title{Bifurcations of Periodic Orbits in the Generalised Nonlinear Schr\"{o}dinger Equation }

\author[1]{Ravindra Bandara\thanks{ravindra@eng.pdn.ac.lk}}
\author[2]{Andrus Giraldo}
\author[3]{Neil G. R. Broderick}
\author[4]{Bernd Krauskopf}
\affil[1]{Department of Engineering Mathematics, Faculty of Engineering, University of Peradeniya 20400, Sri Lanka}
\affil[2]{School of Computational Sciences, Korea Institute for Advanced Studies, Seoul 02455, Korea}
\affil[3]{Department of Physics and Dodd-Walls Centre for Photonic and Quantum Technologies, University of Auckland, Auckland 1010, New Zealand}
\affil[4]{Department of Mathematics and Dodd-Walls Centre for Photonic and Quantum Technologies, University of Auckland, Auckland 1010, New Zealand\\}

\date{}
\maketitle

\begin{abstract}
We focus on the existence and persistence of families of saddle periodic orbits in a four-dimensional Hamiltonian reversible ordinary differential equation derived using a travelling wave ansatz from a generalised nonlinear Schr{\"o}dinger equation (GNLSE) with quartic dispersion. In this way, we are able to characterise different saddle periodic orbits with different signatures that serve as organising centres of homoclinic orbits in the ODE and solitons in the GNLSE. To achieve our objectives, we employ numerical continuation techniques to compute these saddle periodic orbits, and study how they organise themselves as surfaces in phase space that undergo changes as a single parameter is varied. Notably, different surfaces of saddle periodic orbits can interact with each other through bifurcations that can drastically change their overall geometry or even create new surfaces of periodic orbits. Particularly we identify three different bifurcations: symmetry-breaking, period-$k$ multiplying, and saddle-node bifurcations. 
Each bifurcation exhibits a degenerate case, which subsequently gives rise to two bifurcations of the same type that occurs at particular energy levels that vary as a parameter is gradually increased.
 Additionally, we demonstrate how these degenerate bifurcations induce structural changes in the periodic orbits that can support homoclinic orbits by computing sequences of period-$k$ multiplying bifurcations.
\end{abstract}

\section{Introduction}
\label{sec:intro}

Nonlinear generalisations of the Schr{\"o}dinger Equation have served to describe various physical phenomena, including small-amplitude gravity waves in the ocean \cite{Vitanov2013DeepWaterWO}, the dynamics of Bose-Einstein condensates \cite{10.1063/1.4828682}, and the propagation of pulses of light in optical waveguides \cite{GDK}. We focus here on the latter context, where dispersion and Kerr nonlinearity compensate each other to create optical solitons that travel through the waveguide unchanged. The dispersion is generally quadratic to lowest order, but recent studies showed the existence of a novel class of optical solitons, which arise from a delicate balance between the Kerr nonlinearity and quadratic as well as quartic dispersion \cite{GDK,PhysRevA.103.063514, tam2018solitary, tam2019stationary}. They are described by the Generalised Nonlinear Schr{\"o}dinger equation (GNLSE) with quartic dispersion, which takes the form
       \begin{equation}
       \frac{\partial A}{\partial z}=i \gamma |A|^{2} A - i \frac{\beta_{2}}{2}\frac{\partial ^{2} A}{\partial t^{2}}+i\frac{\beta_{4}}{24}\frac{\partial ^{4} A}{\partial t^{4}},
       \label{eq:gnlse}
       \end{equation}
where $A(z, t)$ represents the complex pulse envelope, $\beta_2$ and $\beta_4$ are the coefficients of the quadratic and of the quartic dispersion, respectively, and $\gamma$ denotes the strength of the Kerr nonlinearity. 

The theoretical and numerical study of optical solitons is achieved via the traveling-wave ansatz $A(z,t)=u(t)e^{i\mu z}$, which transforms Eq.~\eqref{eq:gnlse} to the four-dimensional Hamiltonian system \cite{PhysRevA.103.063514, bandara2023generalized, PARKER2021132890}
        \begin{equation}
	\cfrac{d\mathbf{u}}{dt}=f(\mathbf{u};\beta_{2},\beta_{4},\gamma,\mu)=\left( {\begin{array}{cc}u_{2} \\u_{3}\\u_{4}\\ \cfrac{24}{\beta_{4}}\left(\cfrac{\beta_{2}}{2}u_{3}+\mu u_{1}-\gamma u_{1}^3  \right) \\\end{array} } \right).
	\label{eq:4Dode}
	\end{equation}\\
Here $\mathbf{u}=(u_1,u_2,u_3,u_4)=\left(u,\frac{du}{dt},\frac{d^{2}u}{dt^{2}},\frac{d^{3}u}{dt^{3}}\right)$, and the conserved energy function is 	
	 \begin{equation}
 H(\mathbf{u})=u_{2}u_{4}-\frac{1}{2}u_{3}^{2}-\left(\frac{6\beta_{2}u_{2}^{2}-6\gamma u_{1}^{4}+12\mu u_{1}^{2}}{\beta_{4}}\right).
        \label{eq:Hamiltonian}
        \end{equation}
Solitons of the GNLSE \eqref{eq:gnlse} are now given by the $u_1$-component of homoclinic orbits to the origin $\mathbf{0}= (0, 0, 0, 0)$ of system~\eqref{eq:4Dode}, which are trajectories that converge to $\mathbf{0}$ in both forward and backward time \cite{PhysRevA.103.063514,bandara2023generalized, PARKER2021132890}. The overall task is, hence, to find and classify (different types of) homoclinic orbits to $\mathbf{0}$, which requires this equilibrium to be a saddle. Moreover, $H(\mathbf{0}) = 0$ for any value of the parameters of system~\eqref{eq:4Dode}, and this means that any homoclinic orbit to $\mathbf{0}$ lies in the zero-energy level of the Hamiltonian $H$ as well.

Importantly for the structure of its homoclinic orbits, system \eqref{eq:4Dode} features reversibility properties given by the two transformations
	\begin{itemize}
	\item[] $R_{1}: (u_{1},u_{2},u_{3},u_{4}) \rightarrow (u_{1},-u_{2},u_{3},-u_{4})$ and 
	\item[] $R_{2}: (u_{1},u_{2},u_{3},u_{4}) \rightarrow (-u_{1},u_{2},-u_{3},u_{4}),$
\end{itemize}
that is, for any given solution $\mathbf{u}(t)$ of system~\eqref{eq:4Dode} both $R_1(\mathbf{u}(-t))$ and $R_{2}(\mathbf{u}(-t))$ are also solutions \cite{PhysRevA.103.063514, bandara2023generalized}. The points that are (pointwise) invariant under $R_1$ or $R_2$, respectively, form the reversibility sections \cite{champneys1998homoclinic,articleParra,champneys1993hunting}
\begin{itemize}
	\item[] $\Sigma_{1}=\{\mathbf{u}\in \mathbb{R}^{4} : u_{2}=u_{4}=0  \}$ and
	\item[] $\Sigma_{2}=\{\mathbf{u}\in \mathbb{R}^{4} : u_{1}=u_{3}=0  \}$.
\end{itemize}
Any solution of system \eqref{eq:4Dode} that intersects either $\Sigma_{1}$ or $\Sigma_{2}$ is necessarily (setwise) invariant under the corresponding symmetry $R_1$ or $R_2$. We refer to such a solution as $R_1$-symmetric and $R_2$-symmetric, respectively, and as $R^{*}$-symmetric if its is invariant under both $R_1$ and $R_2$; solutions that do not exhibit invariance under $R_1$ or $R_2$ are called non-symmetric \cite{PhysRevA.103.063514,bandara2023generalized}. 

General results for fourth-order, reversible, and Hamiltonian systems \cite{champneys1998homoclinic,devaney1977blue, devaney1976} state that, for any fixed value of the system parameters, periodic orbits come in families that are parameterised locally by the energy $H$. Any homoclinic orbit is the limit of a family of periodic orbits and, moreover, is strucurally stable --- it persists under (sufficiently small) parameter changes. Specifically for system~\eqref{eq:4Dode}, it has been known that a primary homoclinic orbit, which is $R_1$-symmetric, exists throughout the parameter region where $\mathbf{0}$ is a saddle equilibrium, whose eigenvalues may have either real or complex conjugate eigenvalues (with nonzero real parts) \cite{GDK, PhysRevA.103.063514}. The transition between these two cases for a saddle involved in a homoclinic orbit is known as a Belyakov-Devaney (\textbf{BD}) bifurcation, which implies the existence of infinitely many homoclinic orbits of different symmetry types nearby \cite{articleParra}. The \textbf{BD} bifurcation has been identified as the main `generator' in system~\eqref{eq:4Dode} of a plethora of connecting orbits of different types \cite{PhysRevA.103.063514, bandara2023generalized}. Connecting orbits from the equilibrium $\mathbf{0}$ to different periodic orbits in the zero-energy level, known as \emph{EtoP connections}, play a crucial role in the organisation of the different families of homoclinic orbits to $\mathbf{0}$. Each such homoclinic orbit features a certain number of loops near the periodic orbit in question, while converging to $\mathbf{0}$ in forward and backward time \cite{PhysRevA.103.063514}. Moreover, there also exist \emph{PtoP connections} between different periodic orbits, and these can be `combined' with EtoP connections to create additional families of homoclinic orbits to $\mathbf{0}$ that feature different numbers of loops during visits to several periodic orbits; we refer to these solitons as \emph{multi-oscillation solitons} for this reason \cite{bandara2023generalized}. The EtoP and PtoP connections were found to emerge from the point \textbf{BD} in pairs, exists over a range of the quadratic dispersion parameter $\beta_2$ and then vanish at fold bifurcations; here the values of the other parameters are fixed at $\beta_4=-1 < 0$, $\gamma=1 > 0 $ and $\mu=1 > 0$ without loss of generality \cite{PhysRevA.103.063514}. These fold points are accumulation points of the fold bifurcations of associated families of homoclinic orbits to $\mathbf{0}$, which also emerge from the point \textbf{BD} and, hence, exist over very comparable ranges of $\beta_2$. This overall scenario is referred to as \emph{\textbf{BD}-truncated homoclinic snaking} \cite{PhysRevA.103.063514,bandara2023generalized, PARKER2021132890}, and we briefly discuss some of its pertinent features in Sec.~\ref{sec:connections} as the starting point for the work presented here. 

A crucial insight from the above discussion is the following. The existence of any family of homoclinic orbits of stystem \eqref{eq:4Dode} is contingent upon the existence of the corresponding EtoP and PtoP connections. In turn, these connections are entirely determined by the periodic orbits that exist. This emphasises the vital role played by periodic orbits for existence and overall organisation of the homoclinic orbits of system \eqref{eq:4Dode} and, hence, the different kinds of solitons found in the GNLSE \eqref{eq:gnlse}. 

This realisation motivates and necessitates the study of the underlying periodic orbit structure of system \eqref{eq:4Dode}. As was already mentioned, for any fixed value of the parameters, that is, for fixed $\beta_2$, the periodic orbits of a given symmetry type form a surface in $(u_{1},u_{2},u_{3},u_{4})$-space that is locally parameterised by the Hamiltonian energy $H$ as given in~\eqref{eq:Hamiltonian}. Since only the periodic orbits in the zero-energy level may form connections with the saddle equilibrium $\mathbf{0}$, it is important to understand what these surfaces look like when represented in terms of $H$, including how homoclinic orbits to $\mathbf{0}$ form their boundaries. As was pointed out in previous work \cite{PhysRevA.103.063514,bandara2023generalized}, for the intermediate value $\beta_2=0.4$ there are three basic periodic orbits: the $R^*$-symmetric periodic orbit $\Gamma_{*}$ and the pair of $R_1$-symmetric periodic orbits $\Gamma_{1}^{\pm}$, which are mapped to one another by the reversibility $R_2$. When continued, these periodic orbits form the three corresponding basic surfaces $\mathcal{S}_{*}$, $\mathcal{S}_1^{-}$ and $\mathcal{S}_1^{+}$.

We present here a detailed exposition of the structure of the three basic surfaces and how they bifurcate to interact with one another as the parameter $\beta_2$ is increased. To this end, we first consider in Sec.~\ref{sec:surf0p4} the case $\beta_2=0.4$ and present the surfaces $\mathcal{S}_{*}$ and $\mathcal{S}_1^{\pm}$ in the $(u_1, u_2, H)$-space. This allows us to show geometrically how they accumulate in a spiraling manner on the $R_1$-symmetric primary homoclinic orbit and its $R_2$-counterpart, creating infinitely many periodic orbits in the zero-energy level in the process. As we will see, the surfaces $\mathcal{S}_{*}$ and $\mathcal{S}_1^{\pm}$ do not interact with one another for this and lower values of $\beta_2$. However, this changes when $\beta_2$ is increased: these three basic surfaces then connect, which changes the overall geometry and symmetry properties of these surfaces. 

We show in Secs.~\ref{sec:symmbreak} to \ref{sec:nonsymm} that this happens via different bifurcations of periodic orbits, including symmetry-breaking bifurcations, saddle-node bifurcations and period-$k$ multiplying bifurcations of periodic orbits, where the non-trivial Floquet multipliers of the respective periodic orbits have a rational rotation number \cite{GREENE1981468,mackay1982renormalisation}. The bifurcations are encountered at isolated points as a function of the energy $H$. They may coalesce or be created at instances where the respective bifurcation is non-generic in a certain way. Each such degenerate bifurcation occurs at a discrete values of $\beta_2$ and at a specific and positive energy level, and concerns/creates two points of bifurcation either side of the respective level. These degenerate bifurcations are encountered in succession as $\beta_2$ is increased, and they induce changes to the overall geometry of periodic orbits --- including the splitting of certain surfaces and the emergence of new. We study and represent this process by considering the intersection sets of the respective surfaces of periodic orbits with a suitable three-dimensional section. 

Section~\ref{sec:order} then discusses how these changes of geomety with $\beta_2$ affect the periodic orbits within the zero-energy level. Because the degenerate bifurcations occur just above the zero-energy level, one of the created bifurcations rapidly approaches and then crosses the zero-energy level as $\beta_2$ is increased. We identify here the sequences of $\beta_2$-values where this happens --- the finger-print of the degenerate bifurcations in the zero-energy level --- for the $R^*$-symmetric periodic orbit $\Gamma_{*}$. More specifically, we formulate a conjecture on the ordering of the period-$k$ multiplying bifurcations involved and present comprehensive numerical evidence in its support. 

The discovery of sequences of period-$k$ multiplying bifurcations in system \eqref{eq:4Dode} is a central aspect of the work presented here. In fact, we find two distinct types of period-$k$ multiplying bifurcations, each exhibiting a different bifurcation structure as a function of $H$. The first type involves a sequence of period-$k$ multiplying bifurcations that give rise to periodic orbits with $R^*$-symmetry. These periodic orbits exist only for odd values of $k \geq 3$. The second type involves a sequence of period-$k$ multiplying bifurcations that result in the emergence of two pairs of periodic orbits, one with $R_1$ symmetry and the other with $R_2$ symmetry. These periodic orbits exist for both odd and even values of $k \geq 2$. The interplay between these two distinct types of period-$k$ multiplying bifurcations in system \eqref{eq:4Dode} provides a comprehensive explanation of how the overall geometry of the three basic surfaces changes with $\beta_2$. 

The work presented here relies extensively on state-of-the-art computational methods that enable the discovery and continuation of the various periodic orbits, as well as the detection of their bifurcations. Periodic orbits and their associated surfaces are computed via the continuation of solutions of appropriately formulated two-point boundary value problems (2PBVP) \cite{krauskopf2007numerical}. In the context of reversible and Hamiltonian systems, periodic orbits are not isolated in phase space for fixed parameter values. We therefore use the established approach of introducing into system~\eqref{eq:4Dode} an additional parameter that multiplies the gradient of the conserved quantity $H$ \cite{galan2014continuation}. This ensure that solutions of the 2PBVP are isolated, while the new parameter actually remains zero to within machine precision; see also \cite[Sec.~III]{PhysRevA.103.063514}. All computations are performed with the continuation package \textsc{Auto-07p} \cite{doedel2007auto} and its extension \textsc{HomCont} \cite{champneys1996numerical}; the thus obtained periodic orbit data is then visualised within \textsc{Matlab}.

\section{The role of heteroclinic cycles involving periodic orbits}
\label{sec:connections}

\begin{figure}[t!]
   \centering
   \includegraphics{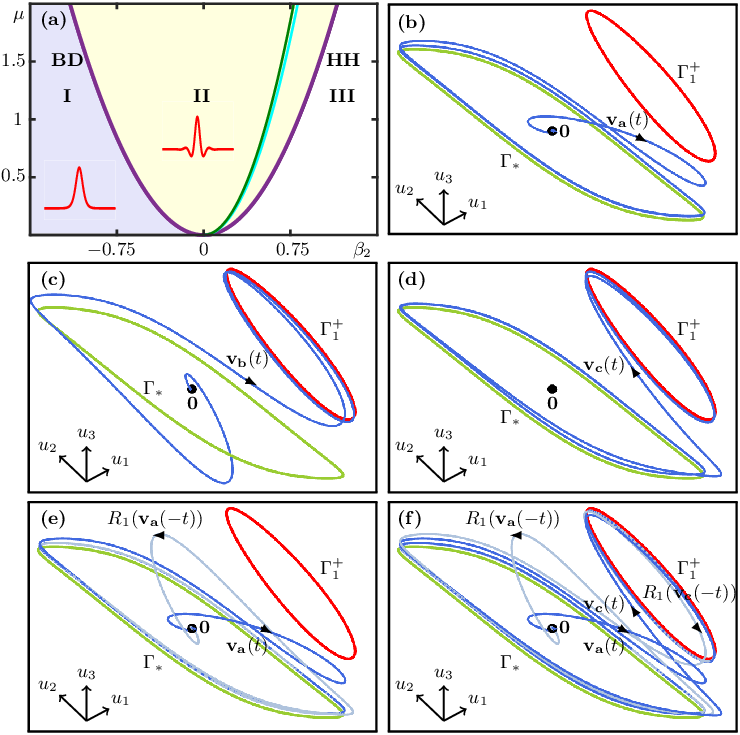}
   \caption{\label{fig:connections} Different types of connecting orbits of system \eqref{eq:4Dode} with $\beta_4=-1$ and $\gamma=1$. Panel~(a) shows the bifurcation diagram in the $(\beta_2, \mu)$-plane with the curves \textbf{BD} and \textbf{HH} delimiting regions I (light purple), II (light yellow) with selected fold curves, and III (white); also shown is the $u_1$-trace of the $R_1$-symmetric primary homoclinic orbit (red curves) in regions I and II, for $(\beta_2, \mu)=(-1,1)$ and $(\beta_2, \mu)=(0.4,1)$, respectively. Panels (b)-(f) are projections onto $(u_1, u_2, u_3)$-space for $\beta_2=0.4$ with $\mu=1$. Panels~(b) and~(c) show EtoP connections (blue curves) from $\mathbf{0}$ (black dot) to $\Gamma_{*}$  (green curve) and $\Gamma_{1}^{+}$ (red curve), respectively, and panel (d) shows a PtoP connection (blue curve) between $\Gamma_{*}$ and $\Gamma_{1}^{+}$. Panel~(e) shows the EtoP cycle between $\mathbf{0}$ and $\Gamma_{*}$, and panel~(f) the heteroclinic cycle from $\mathbf{0}$ to $\Gamma_{1}^{+}$ to $\Gamma_{*}$ and back to $\mathbf{0}$.}
\end{figure} 

\begin{figure}[t!]
   \centering
   \includegraphics{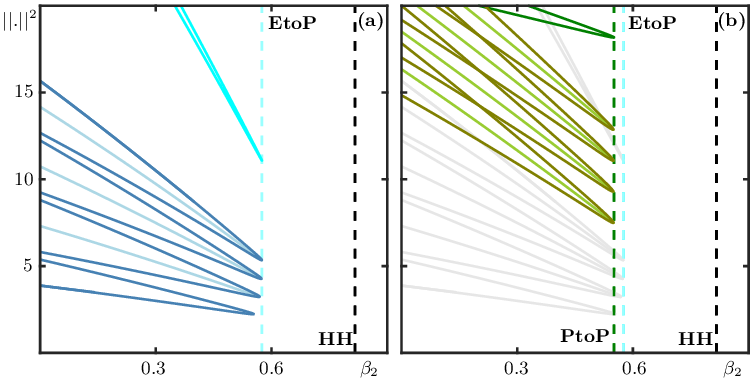}
   \caption{\label{fig:oneparbif} Families of homoclinic orbits associated with the heteroclinic cycles in Fig.~\ref{fig:connections}(e) and (f), respectively. Panel (a) shows the one-parameter bifurcation diagrams in $\beta_2$ with the curve of EtoP connections (cyan) between $\mathbf{0}$ and $\Gamma_{*}$ and curves (dark and light blue) of homoclinic orbits from and back to $\mathbf{0}$ that feature an increasing number of loops near $\Gamma_{*}$. Panel (b) shows additionally the curve of PtoP connections (green) between $\Gamma_{*}$ and $\Gamma_{1}^{+}$  and curves (forest and light green) of homoclinic orbits to $\mathbf{0}$ with three loops near $\Gamma_{*}$, an increasing number of loops near $\Gamma_{1}^{+}$ and again three loops near $\Gamma_{*}$. The vertical dashed lines indicate \textbf{HH} (black) and the fold points of \textbf{EtoP} (cyan) and of \textbf{PtoP} (green).}
\end{figure}

EtoP and PtoP connections involving different types of periodic orbits play a crucial role in generating associated families of homoclinic orbits of system \eqref{eq:4Dode}, as part of the overall phenomenon of \textbf{BD}-truncated homoclinic snaking. We illustrate this in Figs.~\ref{fig:connections} and ~\ref{fig:oneparbif} for the basic periodic orbits $\Gamma_{*}$ and $\Gamma_{1}^{\pm}$; for a comprehensive exposition and further details see \cite{PhysRevA.103.063514, bandara2023generalized}. 

Figure~\ref{fig:connections} shows the bifurcation diagram of system \eqref{eq:4Dode} in the $(\beta_2, \mu)$-plane in panel~(a) and examples of different types of connecting orbits in panels~(b)--(f); here and throughout, we fix $\beta_4=-1$ and $\gamma=1$ without loss of generality \cite{PhysRevA.103.063514}. Figure~\ref{fig:connections}(a) illustrates how the 
$(\beta_2, \mu)$-plane is divided into three regions by the curves \textbf{BD} of Belyakov-Devaney bifurcation and the curve \textbf{HH} of Hamiltonian-Hopf bifurcation. The curves \textbf{BD} and \textbf{HH} meet at the point $(\beta_2, \mu) = (0,0)$ and are the two halves of a parabola; in fact, all curves of bifurcation in the $(\beta_2, \mu)$-plane are semi-parabolas emerging from this point. In regions~I and~II, the point $\mathbf{0}$ is a saddle equilibrium with a pair of $R_1$-symmetric primary homoclinic orbits \cite{PhysRevA.103.063514,tam2018solitary,AKHMEDIEV1994540}, the $u_1$-trace of which for $\beta_2=-1$ and $\beta_2=0.4$ is shown in each of these two regions. At the curve \textbf{BD}, the eigenvalues of $\mathbf{0}$ change from being real in region~I  to being complex conjugate (with nonzero real parts) in region~II, which is the characterising feature of the Belyakov-Devaney bifurcation of a homoclinic orbit \cite{articleParra}. Notice that the $u_1$-trace, that is, the soliton of the GNLSE~\eqref{eq:gnlse}, features oscillating decaying tails in region~II. This is also the region where one finds a plethora of families of secondary homoclinic orbits to $\mathbf{0}$. At the curve \textbf{HH} the eigenvalues of $\mathbf{0}$ become purely complex conjugate (with zero real parts); in region~III this equilibrium is no longer a saddle and, hence, there are no homoclinic orbits to $\mathbf{0}$ in this region. The pair of $R_1$-symmetric primary homoclinic orbits, as well as a pair of $R_2$-symmetric primary homoclinic orbits exist throughout region~II and disappear at the curve \textbf{HH} \cite{PhysRevA.103.063514}. All other, secondary homoclinic solutions emerge at the curve \textbf{BD} bifurcation and disappear, as $\beta_2$ is increased, prior to reaching the \textbf{HH} bifurcation --- namely at curves of fold bifurcations that accumulate on curves of folds of the respective EtoP and PtoP connections organising them. 

Panels~(b) and~(c) of Fig.~\ref{fig:connections} show, in projection onto $(u_1, u_2, u_3)$-space, EtoP connections $\mathbf{v_{a}}(t)$  and $\mathbf{v_{b}}(t)$ of system~\eqref{eq:4Dode} from $\mathbf{0}$ to the basic periodic orbits $\Gamma_{*}$ and $\Gamma_{1}^{+}$, respectively; here $\gamma=1$ (without loss of generality) and we set $\beta_2=0.4$. Notice how these two connecting orbits converge to the equilibrium $\mathbf{0}$ in backward time, while accumulating on the respective periodic orbit in forward time. Similarly, panel~(d) shows a PtoP connection $\mathbf{v_{c}}(t)$ from $\Gamma_{*}$ to $\Gamma_{1}^{+}$ which accumulates on $\Gamma_{*}$ in forward and on $\Gamma_{1}^{+}$ in backward time. Importantly, the respective images of these connecting orbits under the reversibilities $R_1$ and $R_2$ also exist and, hence, provide return connections. Note here that $\mathbf{0}$ and $\Gamma_{*}$ are invariant under both $R_1$ and $R_2$, while $\Gamma_{1}^{+}$ is invariant only under $R_1$ and has the $R_2$-counterpart $\Gamma_{1}^{-}$. 

As Fig.~\ref{fig:connections}(e) and~(f) illustrate, these different connecting orbits can be assembled to more complicated heteroclinic cycles. More specifically, panel~(e) shows the EtoP cycle between $\mathbf{0}$ and $\Gamma_{*}$ formed by $\mathbf{v_{a}}(t)$ and  $R_1(\mathbf{v_{a}}(-t))$, and panel~(f) the heteroclinic cycle consisting of the connections $\mathbf{v_{a}}(t)$ from $\mathbf{0}$ to $\Gamma_{*}$; $\mathbf{v_{c}}(t)$ from $\Gamma_{*}$ to $\Gamma_{1}^{+}$; $R_1(\mathbf{v_{c}}(-t))$ from $\Gamma_{1}^{+}$ back to $\Gamma_{*}$; and $R_1(\mathbf{v_{a}}(-t))$ from $\Gamma_{*}$ back to $\mathbf{0}$. Note that the $R_2$-counterparts of these heteroclinic cycles also exist, which now concerns $\Gamma_{1}^{-}$ rather than $\Gamma_{1}^{+}$.

Figure~\ref{fig:oneparbif} shows how the heteroclinic cycles in Fig.~\ref{fig:connections}(e) and (f) give rise to infinite families of homoclinc orbits to $\mathbf{0}$. The EtoP connection between $\mathbf{0}$ and $\Gamma_{*}$ from Fig.~\ref{fig:connections}(b) can be continued in the parameter $\beta_2$ as a single curve with a fold $\beta_2\approx0.5752$. It is shown in the one-paramter diagram in Fig.~\ref{fig:oneparbif}(a) together with curves of associated homoclinc orbits to $\mathbf{0}$ that feature more and more loops near $\Gamma_{*}$. Here all computed connecting orbits are represented by their $L_2$-norm $||.||^2$, and dashed vertical lines indicate the \textbf{HH} bifurcation and the fold \textbf{EtoP} of the EtoP connection. Note that the curves of homoclinc orbits in panel~(a) all have folds that converge quickly to \textbf{EtoP} as the number of loops near $\Gamma_{*}$ increases, which is represented by higher values of their $L_2$-norm. For images of representative $u_1$-trace for this specific family of homoclinic orbits, which are intricately associated with the solitons found in the GNLSE, we refer the reader to \cite[Fig.~4(c)--(f)]{bandara2023generalized}.

The one-parameter bifurcation diagram in $\beta_2$ in {Fig.~\ref{fig:oneparbif}(b) shows the continuation of the PtoP connection between $\Gamma_{*}$ and $\Gamma_{1}^{+}$ from Fig.~\ref{fig:connections}(d); it is again a single curve with a fold at $\beta_2\approx0.551$, which is marked by the dashed vertical line labeled \textbf{PtoP}; for ease of comparison, we show \textbf{EtoP} and the curves from panel~(a) in the background as well. Also shown in panel~(b) are curves of associated homoclinic orbits from $\mathbf{0}$ that feature three loops near $\Gamma_{*}$, an increasing number of loops near $\Gamma_{1}^{+}$ and another three loops near $\Gamma_{*}$, before returning back to $\mathbf{0}$. These families of homoclinic orbits require that both the EtoP and the PtoP exist, and their respective curves all have folds whose $\beta_2$-values are indistinguishable from that of the fold \textbf{PtoP}, which has the lower $\beta_2$-value compared to \textbf{EtoP}; as Fig.~\ref{fig:connections}(a) shows, the limiting loci \textbf{EtoP} and \textbf{PtoP} form semi-parabolas in the $(\beta_2, \mu)$-plane of system~\eqref{eq:4Dode}. Images of $u_1$-traces of these families of homoclinic orbits can be found in \cite[Fig.~12(a)--(f)]{bandara2023generalized}. We remark that the lighter curves in panels~(a) and~(b) represent symmetry-broken homoclinic orbit that bifurcate from the respective fold points in symmetry-breaking bifurcations; they are associated with non-symmetric heteroclinic cycles involving both simultaneously existing EtoP and PtoP connections, which are not related by symmetry. 

In summary, Figs.~\ref{fig:connections} and~\ref{fig:oneparbif} give an impression of how connecting orbits involving the basic periodic orbits $\Gamma_{*}$ and $\Gamma_{1}^{\pm}$ can be combined to give rise to associated infinite families of homoclinic orbits to $\mathbf{0}$ of system~\eqref{eq:4Dode}. For a comprehensive presentation of the different types of connecting and homoclinic orbits and their symmetry properties we refer the reader to \cite{PhysRevA.103.063514, bandara2023generalized}. The overall message is that, in the same way, any other periodic orbit in the zero-energy level can be included in ever more complex heteroclinic cycles from and back to equilibrium $\mathbf{0}$ --- to create a never ending menagerie of families of homoclinic orbits \cite{Lohse2016BoundaryCF, PARKER2022368} and, hence, a huge variety of solitons of the GNLSE~\eqref{eq:gnlse}. This is why we next study the structure of periodic orbits of system~\eqref{eq:4Dode} and how it changes with the parameter $\beta_2$.

\section{Geometry of the three basic surfaces for $\beta_2=0.4$}
\label{sec:surf0p4}

Continuation of the $R^{*}$-symmetric periodic orbit $\Gamma_{*}$ and the $R_{1}$-symmetric periodic orbit $\Gamma_{1}^{+}$, with the continuation package \textsc{Auto-07p} \cite{doedel2007auto} and the 2PBVP setup from \cite[Sec.~III]{PhysRevA.103.063514}, allows us to find the associated one-parameter families of periodic orbits for fixed $\beta_2=0.4$. By rendering them as surfaces $\mathcal{S}_{*}$ and $\mathcal{S}_{1}^{+}$ (and $\mathcal{S}_{1}^{-}$ as its image under $R_2$) we are able to show why and how they contain infinitely many further periodic orbits in the zero-energy level as they accumulate on the pair of $R_{1}$-symmetric primary homoclinic orbits. We first consider the $R^{*}$-symmetric surface $\mathcal{S}_{*}$, then the $R_1$-symmetric surface $\mathcal{S}_{1}^{+}$, and finally show all three basic surfaces together.

\subsection{The $R_*$-symmetric surface $\mathcal{S}_{*}$}
\label{sec:surfS*}

\begin{figure}[t!]
   \centering
   \includegraphics{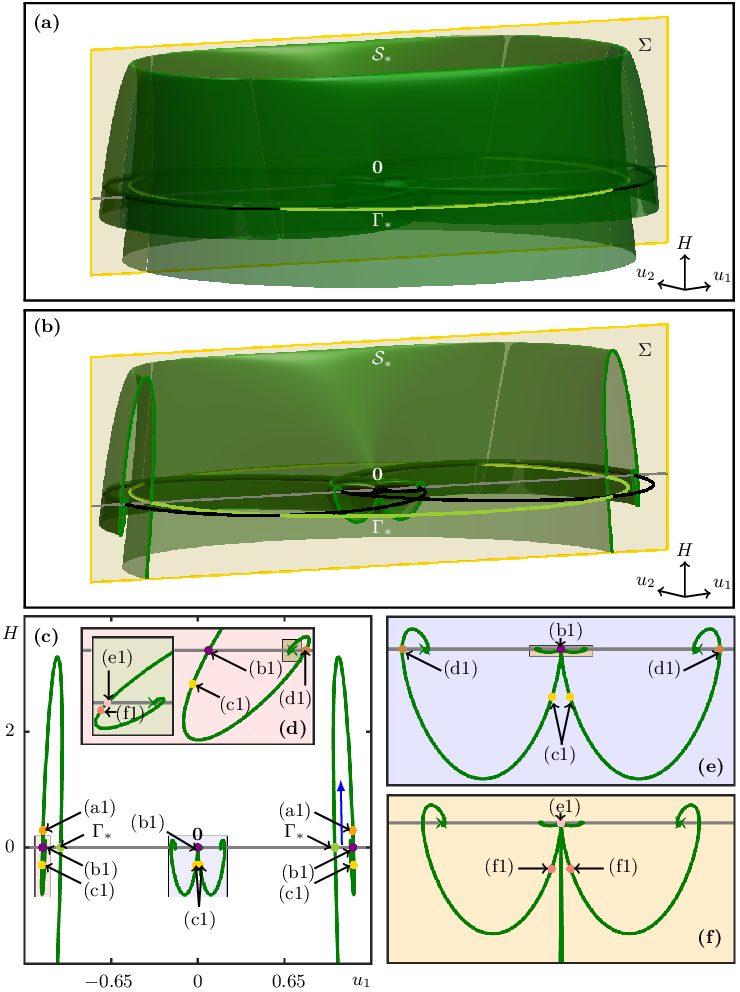}
   \caption{\label{fig:S*} The surface $\mathcal{S}_{*}$ for $\beta_2=0.4$, shown in $(u_1, u_2, H)$-space (a) and as a cutaway view (b), with $\mathbf{0}$ (black dot), $\Gamma_{*}$ (light green curve), the $R_1$-symmetric primary homoclinic orbits (black curves), and the section $\Sigma$ defined by $u_2=0$ (beige plane) with the zero-energy level (gray line). Panel (c) and the enlargements~(d)--(f) show $\mathcal{S}_{*} \cap \Sigma$ in the $(u_1,H)$-plane with the primary homoclinic orbits (green crosses) and periodic orbits (coloured dots, labeled) shown in Fig.~\ref{fig:tanS*}; period increases in the direction of the blue arrow.}
\end{figure} 

Figure \ref{fig:S*} shows the surface $\mathcal{S}_{*}$ of $R^*$-symmetric periodic orbits of system \eqref{eq:4Dode} for $\beta_2=0.4$ in different ways. Panel~(a) is an overall view in projection onto $(u_{1},u_{2},H)$-space, where we also show the pair of $R_1$-symmetric primary homoclinic orbits, the basic periodic orbit $\Gamma_{*}$, and the section $\Sigma$ defined by $u_2=0$ and appearing as a plane in this projection. Notice that the surface $\mathcal{S}_{*}$ is indeed symmetric with respect to both $R_1$ and $R_2$, that is under changing the sign of either $u_2$ or $u_1$. Figure \ref{fig:S*}(b) shows the same objects but now, in a cutaway view, we only render the part of $\mathcal{S}_{*}$ where $u_2 \geq 0$; its image under $R_1$ is the part of $\mathcal{S}_{*}$ that is not shown. The advantage of this representation in $(u_{1},u_{2},H)$-space is that each periodic orbit is seen to lie in a specific energy level. The surface $\mathcal{S}_{*}$ has a global maximum when it reaches a periodic orbit with $H(\mathbf{u})\approx 3.3$. However, it does not have a global minimum in $H$, and our numerical continuation results strongly suggest that the surface $\mathcal{S}_{*}$ extends to any negative value of $H$. Note that $\Gamma_{*}$ is the `first' periodic orbit of $\mathcal{S}_{*}$ in the zero-energy level when increasing $H$ from large negative values, and it only makes single loop around the point $\mathbf{0}$. 

As the surface $\mathcal{S}_{*}$ accumulates on the pair of $R_1$-symmetric primary homoclinic orbits, the period increases and one encounters further periodic orbits in the zero-energy level. This is illustrated in Fig.~\ref{fig:S*}(c) by showing the intersections set $\mathcal{S}_{*} \cap \Sigma$  
in projection onto the $(u_1,H)$-plane. Notice that the intersection set $\mathcal{S}_{*} \cap \Sigma$ has three components: two outer spirals, which are associated with the `far' intersection points of the two $R_1$-symmetric primary homoclinic orbits (and each others image under $R_2$), and a central component near the point $\mathbf{0}$. Successive enlargements of the left spiral in panel~(d) and its inset illustrate the spiraling nature of $\mathcal{S}_{*} \cap \Sigma$ as it approaches the left-most intersection point of the basic homoclinic orbits; hence, there is a sequence of $R^*$-symmetric periodic orbits of increasing period in the zero-energy level. As the successive enlargements of the central region in Fig.~\ref{fig:S*}(d) and~(e), show, the set $\mathcal{S}_{*} \cap \Sigma$ develops further branches near $\mathbf{0}$ in the process, which spiral into subsequent intersection points of the pair of basic homoclinic orbits.

In order to explain how the sequence of periodic orbits in the zero-energy level and this geometry of $\mathcal{S}_{*} \cap \Sigma$ arise, Fig.~\ref{fig:tanS*} shows in $(u_1,u_2,u_4)$-space a number of selected periodic orbits on $\mathcal{S}_{*}$, namely those that are highlighted and accordingly labeled in Fig.~\ref{fig:S*}(c)--(f). Panels~(a1)--(f1) of Fig.~\ref{fig:tanS*} show how the periodic orbit on $\Sigma$ changes as its period is increased; also shown is the section $\Sigma$, which appears again as a plane, and the $R_1$-reversibility section $\Sigma_1 \subset \Sigma$, which appears as a line in this projection. Starting from $\Gamma_{*}$ --- the first periodic orbit in the zero-energy level --- and increasing the period we find first the single-loop periodic orbit in panel~(a1) with $H = 0.3$; it has pair of intersection points $P_1$ and $P_2$ with $\Sigma_1$ and also a pair of small loops, one either side of $\Sigma$. Indeed, the fact that these objects come in pairs is due to the $R^*$-symmetry, which manifests itself as a rotation around the $u_4$-axis in Fig.~\ref{fig:tanS*}. When the period increased further the small loops grow and the energy decreases again until the two loops become tangent to the section $\Sigma$. This situation is shown in panel~(b1) and it occurs exactly at $H=0$, that is, for the second periodic in the zero-energy level; see also the enlargement in panel~(b2) and notice that this pair of tangencies does not occur on the $R_1$-reversibility section $\Sigma_1$. For larger values of the period, as is shown for $H = -0.3$ in panels~(c1) and~(c2), this leads to pairs of two additional intersection points with the section $\Sigma$, which we label $P_3$ to $P_6$. Notice also that these loops now `encircle' the equilibrium $\mathbf{0}$. Observe that the (pairs of) intersection points $P_3$/$P_5$ and $P_4$/$P_6$ generate the largest additional branches of $\mathcal{S}_{*} \cap \Sigma$ in Fig.~\ref{fig:S*}(c) and~(d), which hence exist after the tangency where $H=0$ for the second time.

\begin{figure}[t!]
   \centering
   \includegraphics{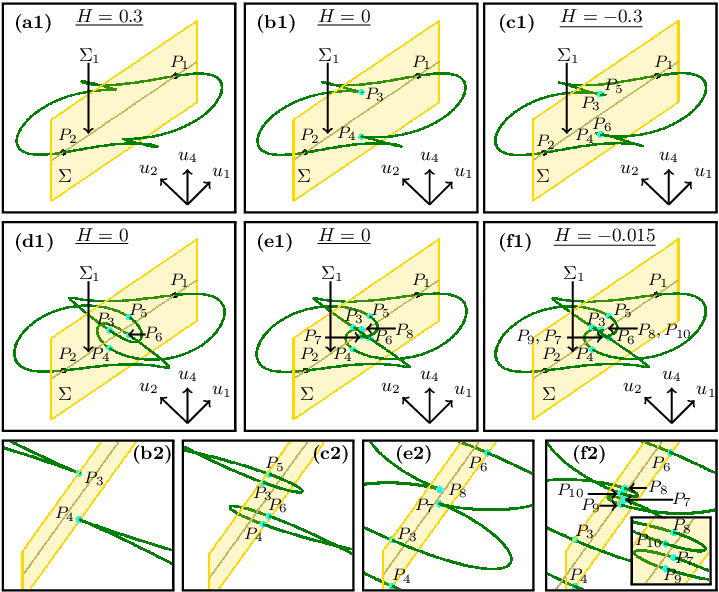}
   \caption{\label{fig:tanS*} Selected periodic orbits (green curves) on the surface $\mathcal{S}_{*}$, as labeled in Fig.~\ref{fig:S*}(c)--(f), shown in $(u_1,u_2,u_4)$-space with the section $\Sigma$; also shown is the $R_1$-reversibility section $\Sigma_1 \subset \Sigma$ that contains two intersection points (black dots). Panels~(a1)--(a6) are overall views (with the respective value of $H$ shown), and panels (b2), (c2), (e2) and (f2) are enlargements at/near quadratic tangencies of the periodic orbit with $\Sigma$, which generate additional intersection points (cyan dots).} 
\end{figure} 

The pair of extra loops grows with the period but only up to a certain size when $H=0$ for the third time as in Fig.~\ref{fig:tanS*}(d1). Subsequently, when the period is increased further, there are again smaller loops developing, namely on the pairs of previous loops, and these become tangent to $\Sigma$ off $\Sigma_1$ as well when $H=0$ for the fourth time; this hard to see in panel~(e1), but clearly visible in the enlargement panel~(e2). As panels~(f1) and~(f2) show, two further pairs $P_7$/$P_9$ and $P_8$/$P_{10}$ of intersection points of the periodic orbit with $\Sigma$ therefore emerge for even larger period. Indeed, these four extra intersection points correspond to the second largest additional branches of $\mathcal{S}_{*} \cap \Sigma$ that are the focus of in Fig.~\ref{fig:S*}(f). Moreover, the new loops now also encircle $\mathbf{0}$. We remark that our numerical evidence from the continuation of the periodic orbits on $\mathcal{S}_{*}$ clearly suggest that their tangencies with $\Sigma$ are quadratic, that is, generic; the fact that in Fig.~\ref{fig:tanS*}(b2) and~(e2) they have the appearance of \say{kinks} \cite{Kent_1992} is entirely due to the shown projection. 

The emergence of additional, smaller loops on previous loops and their subsequent tangency with $\Sigma$ repeats ad infinitum with every full rotation around (the main) spirals of $\mathcal{S}_{*} \cap \Sigma$ in $(u_1, u_2, H)$-space. In fact, this represents the process of how the family of periodic orbits on $\mathcal{S}_{*}$ converges, for increasing period, to the pair of $R_1$-symmetric primary homoclinic orbits and, in particular, also accumulates on the equilibrium $\mathbf{0}$ (with complex conjugate eigenvalues in region~II). As we have seen, this not only generates an infinite sequence of periodic orbits in the zero-energy level, but also leads to the emergence of ever more branches of the central part of $\mathcal{S}_{*} \cap \Sigma$. These branches correspond to periodic orbits with different numbers of loops around $\mathbf{0}$; note here that these numbers are odd, because $\Gamma_{*}$ makes a single loop around $\mathbf{0}$ and any additional loops around $\mathbf{0}$ are generated at the tangencies with $\Sigma$ in pairs due to the $R^{*}$-symmetry of the surface $\mathcal{S}_{*}$.

\subsection{The $R_1$-symmetric surface $\mathcal{S}_{1}^{+}$}
\label{sec:S+}

\begin{figure}[t!]
   \centering
   \includegraphics{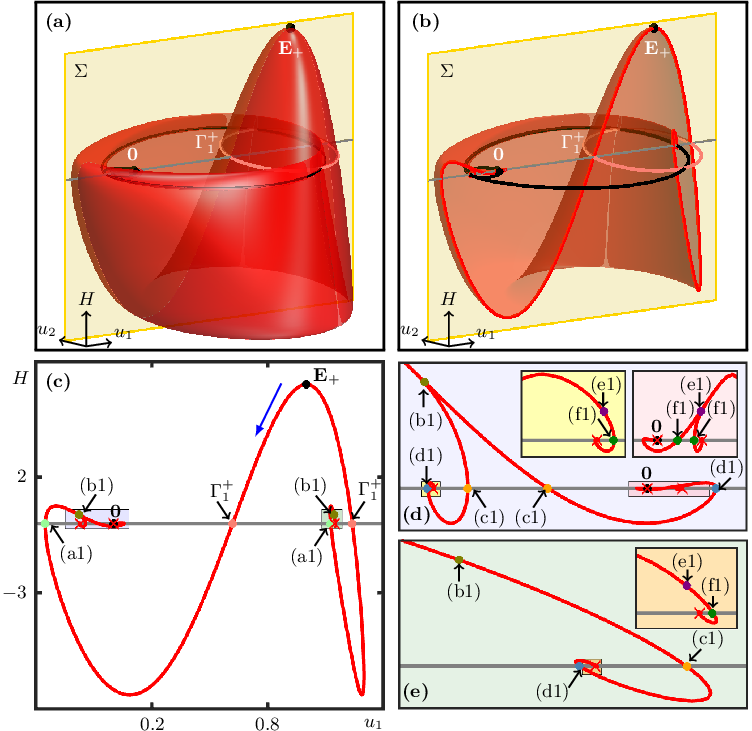}
    \caption{\label{fig:S1p} The surface $\mathcal{S}_1^{+}$ for $\beta_2=0.4$, shown in $(u_1, u_2, H)$-space (a) and as a cutaway view (b), with $\mathbf{0}$ and $\mathbf{E}_{+}$ (black dots), the right-most $R_1$-symmetric primary homoclinic orbits (black curve), the periodic orbit $\Gamma_{1}^+$ (light red curve), and the section $\Sigma$ defined by $u_2=0$ (beige plane) with the zero-energy level (gray line). Panel (c) and successive enlargements~(d) and~(f) show $\mathcal{S}_{1}^+ \cap \Sigma$ in projection onto the $(u_1,H)$-plane with the intersection points of the $R_1$-symmetric primary homoclinic orbit (green crosses) and selected periodic orbits (coloured dots, labeled) shown in Fig.~\ref{fig:tanS1p}; period increases from $\mathbf{E}_{+}$ in the direction of the blue arrow.}  
\end{figure} 

The surface $\mathcal{S}_{1}^{+}$ of $R_1$-symmetric periodic orbits as been found by continuation from the periodic orbit $\Gamma_{1}^{+}$, and it is shown in Figs.~\ref{fig:S1p} in the style of Figs.~\ref{fig:S*}. Note from the two views in $(u_1, u_2, H)$-space in panels~(a) and~(b) of Fig.~\ref{fig:S1p} that $\mathcal{S}_{1}^{+}$ is not invariant under $R_2$, lies mostly in the region of positive $u_1$, and accumulates on the respective $R_1$-symmetric primary homoclinic orbit, which is also shown. The surface $\mathcal{S}_1^{+}$ has a global maximum $H(\mathbf{E}_+)$ in $H$ at the equilibrium $\mathbf{E_+}$, from which the family of periodic orbits emerges. The first periodic orbit in the zero-energy level, as the period is increased, is the basic periodic orbit $\Gamma_{1}^{+}$. When the period is increased further, a global minimum in $H$ is reached in the form of a periodic orbit with $H(\mathbf{u}) \approx -7.4$. Thus, the (Hamiltonian) energy of this family of $R_1$-symmetric periodic orbits is bounded between these two values. This is indeed also the case for the $R_2$-counterpart $\mathcal{S}_{1}^{-}$ (not shown in Fig.~\ref{fig:S1p}) in the region of mostly negative $u_1$, where the maximum of $H$ occurs at the corresponding equilibrium $\bf{E_-}$ and the minimum at the corresponding periodic orbit.

Fig.~\ref{fig:S1p}(c) and the further enlargements in panels~(d) and~(e) show the projection onto the $(u_1,H)$-plane of the intersection sets of these objects with the section $\Sigma$ with $u_2=0$. The set $\mathcal{S}_{1}^{+} \cap \Sigma$ has two main branches, one to the left and one to the right of the point $\mathbf{E_+} \in \Sigma$. The periodic orbit $\Gamma_{1}^{+}$ has exactly one intersection point on each of these two main branches, and this is also the case for the next periodic orbit with $H=0$. However, $\mathcal{S}_{1}^{+} \cap \Sigma$ develops additional branches at specific point on the left main branch as the surface $\mathcal{S}_{1}^{+}$ accumulates on the $R_1$-symmetric primary homoclinic orbit when the period is increased further. All these branches sprial into intersection points of the $R_1$-symmetric primary homoclinic orbit with $\Sigma$; see Fig.~\ref{fig:S1p}(d) and~(e).

\begin{figure}[t!]
   \centering
   \includegraphics{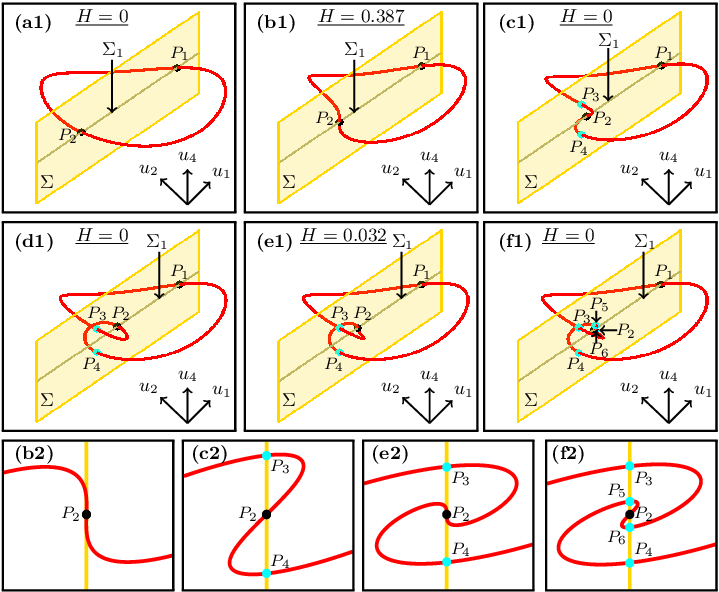}
    \caption{\label{fig:tanS1p} Selected periodic orbits (green curves) on the surface $\mathcal{S}_{1}^{+}$, as labeled in Fig.~\ref{fig:S*}(c)--(f), shown in $(u_1,u_2,u_4)$-space with the section $\Sigma$; also shown is the $R_1$-reversibility section $\Sigma_1 \subset \Sigma$ that contains two intersection points (black dots). Panels~(a1)--(a6) are overall views (with the respective value of $H$ shown), and panels (b2), (c2), (e2) and (f2) are enlargements at/near cubic tangencies of the periodic orbit with $\Sigma$, which generate additional intersection points (cyan dots).}
\end{figure} 

Figure~\ref{fig:tanS1p} illustrates in the $(u_1,u_2,u_4)$-space that new branches emerge at a sequence of cubic tangencies of the periodic orbit with the section $\Sigma$ as the period increases; here the panels~(a1)--(f1), which also show the respective values of $H$, correspond to the intersection sets of the periodic orbits that are highlighted and labeled in Fig.~\ref{fig:S1p}(c)--(e). The enlargement panels~(b2),~(c2) and~(e2),~(f2) of  Fig~\ref{fig:tanS1p} show the projection onto the $(u_2,u_4)$-plane, so that $\Sigma$ appears as a vertical line. Fig~\ref{fig:tanS1p}~(a1) shows the next periodic orbit in the zero-energy level, as we increase the period from $\Gamma^{+}$. It is still a single-loop periodic orbit that intersects $\Sigma$ in two points $P_1$ and $P_2$, which actually lie in the reversibility section $\Sigma_1 \subset \Sigma$. Past the second local maximum of $\mathcal{S}_{1}^{+} \cap \Sigma$ one encounters a cubic tangency of the periodic orbit on $\mathcal{S}_{1}^{+}$ with $\Sigma$ at the point $P_2 \in \Sigma_1$; it is illustrated in Fig~\ref{fig:tanS1p}~(b1) and the enlargement~(b2). When the period in increased, the cubic tangency gives rise to two additional intersetion points $P_3$ and $P_4$, which do not lie in the reversibility section $\Sigma_1$, but either side of it in $\Sigma$. Panels~(c1) and~(c2) show this situation for the periodic with $H=0$. Comparison with Fig~\ref{fig:S1p}(d) shows that this first cubic tangency indeed gives rise to the first additional branch of the intersection set $\mathcal{S}_{1}^{+} \cap \Sigma$; note here, that $P_3$ and $P_4$ have the same values of $u_1$ and $H$ due to the $R_1$-invariance of the surface $\mathcal{S}_{1}^{+}$, so that they appear as one and the same branch in the shown projection of $\Sigma$ onto the $(u_1,H)$-plane. The next periodic orbit in the zero-energy level, shown in Fig~\ref{fig:tanS1p}~(d1), is therefore clearly a two-loop orbit. As the periodic is increased further, we find the second cubic tangency illustrated in panels~(e1) and~(e2), which again happens at $P_2$ and lead to two further intersection points $P_5$ and $P_6$ either side of $\Sigma_1$; the next such orbit in the zero energy surface is shown in panels~(f1) and~(f2). 

The respective periodic orbits in Fig~\ref{fig:tanS1p}(b) and~(e) are tangent to $\Sigma$ at $P_2$ in the direction of $u_4$ and, hence, transverse to the reversibility section $\Sigma_1 \subset \Sigma$; this is numerical evidence that they are generic cubic tangencies. At $P_1 \in \Sigma_1$ and at all further intersection points $P_i \not\in \Sigma_1$, on the other hand, the respective periodic orbits intersect $\Sigma$ transversely throughout. Moreover, our continuation of the periodic orbits from $\Gamma_{1}^+$ clearly suggests that there is an infinite sequence of generic cubic tangencies at $P_2$ that increases the number of loops the periodic orbits on $\mathcal{S}_{1}^{+}$ as they accumulate on the $R_1$-symmetric primary homoclinic orbit. This explains to the observed emergence of additional branches in Fig~\ref{fig:S1p}(c) and~(d). 

\subsection{Geometry of $\mathcal{S}_{*}$ and $\mathcal{S}_{1}^{\pm}$ relative to one another}
\label{sec:geomSall}

\begin{figure}[t!]
   \centering
   \includegraphics{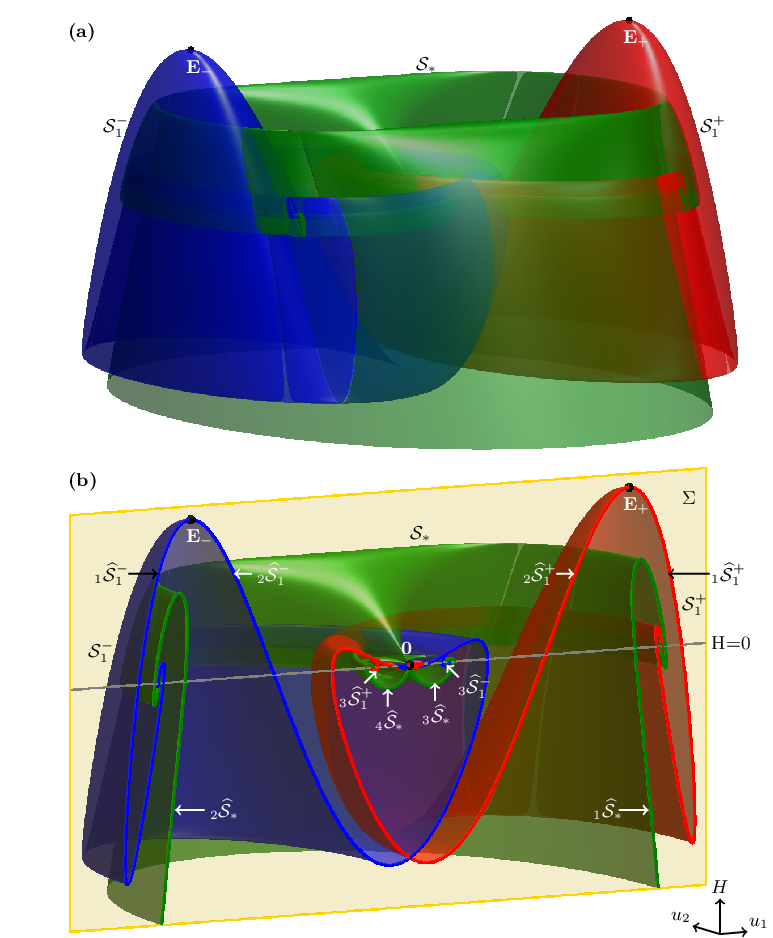}
    \caption{\label{fig:Sboth} The surfaces $\mathcal{S}_{*}$ (green), $\mathcal{S}_{1}^{+}$ (red) and $\mathcal{S}_{1}^{-}$ (blue) for $\beta_2=0.4$, shown in $(u_1, u_2, H)$-space with $\mathbf{0}$ and $\mathbf{E}_{\pm}$ (black dots). Panel~(a) shows the entire surfaces, and panel~(b) is a cutaway view with the section $\Sigma$ defined by $u_2=0$ (beige plane) and the zero-energy level (gray line).}
\end{figure} 

Figure~\ref{fig:Sboth} shows the three basic surfaces $\mathcal{S}_{*}$ and $\mathcal{S}_{1}^{\pm}$ for $\beta_2=0.4$ together in $(u_1,u_2,H)$-space; also shown are the equilibria $\mathbf{0}$ and $\mathbf{E}_{\pm}$. Note that $\mathcal{S}_{1}^{-}$ is obtained as the $R_2$ counterpart of $\mathcal{S}_{1}^{+}$ and, hence, has the same properties. These three surfaces do not intersect, but all accummulate on the $R_1$-symmetric primary homoclinic orbits. This is not obvious from the overall view in Fig.~\ref{fig:Sboth}(a), since there are visible intersections. However, these are due to the projection onto $(u_1,u_2,H)$-space, as the cutaway view in panel~(b) shows. More specifically, the three intersection sets in the section $\Sigma$, which is represented here again by the $(u_1, H)$-plane, do indeed not intersect as they all spiral into the (infinitely many) intersection points of the $R_1$-symmetric primary momoclinic orbits with $\Sigma$; compare with Figures~\ref{fig:S*}(b) and~\ref{fig:S1p}(b).

As we will see in the next sections, the three basic surfaces interact as the parameter $\beta_2$ is increased. Since this is best studied by considering specific branches of the corresponding intersection sets with $\Sigma$, we distinguish and label them in Fig.~\ref{fig:Sboth}(b). Here and throughout, the intersection set of any surface is denoted with a hat symbol, and we use left subscripts to enumerate and distinguish its different branches. In this way, the two main outer spirals of $\mathcal{\widehat{S}}_{*}$ in Fig.~\ref{fig:Sboth}(b) are ${}_{1}\mathcal{\widehat{S}}_{*}$ and  ${}_{2}\mathcal{\widehat{S}}_{*}$, and further pairs of branches are labeled ${}_{3}\mathcal{\widehat{S}}_{*}$, ${}_{4}\mathcal{\widehat{S}}_{*}$ and so on; compare with Figure~\ref{fig:S*}(c)--(e). The two main branches of $\mathcal{\widehat{S}}_{1}^{\pm}$ that emerge from $\mathbf{E}_{\pm}$ are labeled ${}_{1}\mathcal{\widehat{S}}_{1}^{\pm}$ and ${}_{2}\mathcal{\widehat{S}}_{1}^{\pm}$ in Fig.~\ref{fig:Sboth}(b), and the further pairs branches that emerge in pairs are enumerated accordingly; compare with Figure~\ref{fig:S1p}(c)--(e). We remark that the branches ${}_{3}\mathcal{\widehat{S}}_{*}$, ${}_{4}\mathcal{\widehat{S}}_{*}$ and ${}_{5}\mathcal{\widehat{S}}_{*}$, ${}_{6}\mathcal{\widehat{S}}_{*}$, as well as the branches ${}_{3}\mathcal{\widehat{S}}_{*}$ and ${}_{4}\mathcal{\widehat{S}}_{*}$, coincide in the shown projection of $\Sigma$ onto the $(u_1, H)$-plane.

\section{Symmetry-breaking bifurcation of $R^{*}$-symmetric periodic orbits}
\label{sec:symmbreak}

\begin{figure}[t!]
   \centering
   \includegraphics[width=11.8cm]{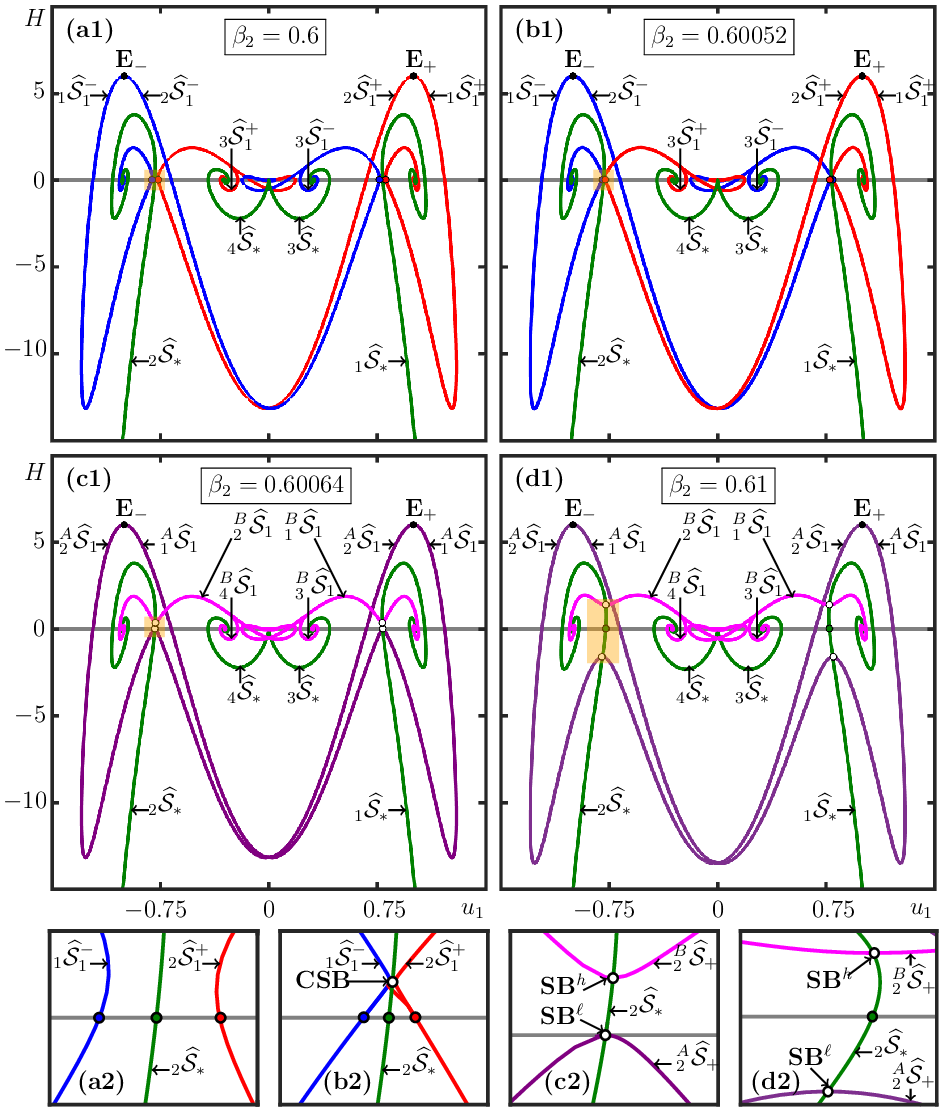}
    \caption{\label{fig:Strans} Creation of symmtry-breaking bifurcations. Panels~(a1)--(d1) show for the stated $\beta_2$-values the intersection sets of the respective surfaces with $\Sigma$ in projection onto the $(u_1,H)$-plane with the zero-enery level (gray line), and panels~(a2)--(d2) are enlargements of the highlighted regions. The branches in panels~(a) are those from Fig.~\ref{fig:Sboth}(b), and they meet at the point $\mathbf{CSB}$ (white dot) in panels~(b). Panels~(c) and~(d) show how the intersection sets of the new surfaces ${}^{A}\mathcal{S}_{*}$ (purple curves) and ${}^{B}\mathcal{S}_{*}$ (magenta curves), which intersect $\mathcal{\widehat{S}}_{*}$ at the points labeled $\mathbf{SB^{\textit{h}}}$ and $\mathbf{SB^{\ell}}$ (white dots).  Also shown are $\mathbf{0}$ and $\mathbf{E}_{+}$ (black dots), while red, blue and green dots, respectively, indicate periodic orbit of $\mathcal{S}_{1}^{+}$, $\mathcal{S}_{1}^{-}$ and $\mathcal{S}_{*}$ in the zero-energy level.} 
\end{figure} 

The first qualitative change in the structure of periodic orbits, as $\beta_2$ is increased, occurs at a degenerate symmetry breaking bifurcation, which we refer to as $\mathbf{CSB}$. At this parameter point, the surfaces $\mathcal{S}_{1}^{+}$ and $\mathcal{S}_{1}^{-}$ connect with the surface $\mathcal{S}_{*}$, namely at a degenerate symmetry-breaking periodic orbit. To keep notation managable, we employ the symbol $\mathbf{CSB}$ to represent both the bifurcation itself as well as this special periodic orbit in phase space. We adopt this convention also for any further bifurcations and their emerging periodic orbits, starting with the symmetry breaking bifurcations/periodic orbits $\mathbf{SB^{\textit{h}}}$ and $\mathbf{SB^{\ell}}$ emerging from $\mathbf{CSB}$.

This transition and its consequences are illustrated in Fig.~\ref{fig:Strans} by showing the intersection curves of the corresponding surfaces in projection onto the $(u_1,H)$-plane. The intersection sets for $\beta_2=0.6$ in panels~(a) are those of the three basic surfaces $\mathcal{S}_{*}$ and $\mathcal{S}_{1}^{\pm}$ for $\beta_2=0.4$, and their branches are labeled as introduced in Fig.~\ref{fig:Sboth}(b). As is illustrated in Fig.~\ref{fig:Strans}(b), for $\beta_2 \approx 0.60052$ the branches ${}_{1}\mathcal{\widehat{S}}_{1}^{-}$, ${}_{2}\mathcal{\widehat{S}}_{1}^{+}$ and ${}_{2}\mathcal{\widehat{S}}_{*}$ meet in the energy level with $H \approx 0.1652$ at a common point labeled $\mathbf{CSB}$, and so do ${}_{2}\mathcal{\widehat{S}}_{1}^{-}$, ${}_{2}\mathcal{\widehat{S}}_{1}^{+}$ and ${}_{1}\mathcal{\widehat{S}}_{*}$ due to $R_1$-symmetry. When $\beta_2$ is increased, the respective intersection sets of $\mathcal{S}_{1}^{-}$ and $\mathcal{{S}}_{1}^{+}$ connect differenly in $\Sigma$ to form two new surfaces ${}^{A}\mathcal{{S}}_{*}$ and ${}^{B}\mathcal{{S}}_{*}$, which meet the surface $\mathcal{{S}}_{*}$ along two symmetry-breaking periodic orbits $\mathbf{SB^{\textit{h}}}$ and $\mathbf{SB^{\ell}}$, with higher and lower energy than $H \approx 0.1652$, respectively. This new structure of the periodic orbits is illustrated in panels~(c) and~(d) of Fig.~\ref{fig:Strans} at the level of the corresponding intersection sets. Here we again enumerate and distinguish by left subscripts the different branches of the intersections sets ${}^{A}\mathcal{\widehat{S}}_{1}$ and ${}^{B}\mathcal{\widehat{S}}_{1}$. We remark that the right subscript in this notation indicates that these surfaces consist of only $R_1$-symmetric periodic orbits. It is important to keep in mind, however, that these surfaces are invariant (as surfaces) under both $R_1$ and $R_2$ as they consist of $R_1$-symmetric periodic orbits and their corresponding $R_2$-counterparts. Indeed, this is the first time that we encounter a difference between the symmetry properties of periodic orbits of a given family and those of the corresponding surface they generate.

For fixed $\beta_2$, as in Fig.~\ref{fig:Strans}(c2) and~(d2), there is a generic symmetry-breaking bifurcation at both $\mathbf{SB^{\textit{h}}}$ and $\mathbf{SB^{\ell}}$ of $R^*$-symmetric periodic orbits that gives rise of a pair of $R_2$-symmetry broken (that is, only $R_1$-symmetric) periodic orbits. Importantly, the bifurcation parameter here is the Hamiltonian energy $H$ (or alternatively the period of the periodic orbits). In particular, the existence of the symmetry-breaking periodic orbits $\mathbf{SB^{\textit{h}}}$ and $\mathbf{SB^{\ell}}$ is a structurally stable property under variation of system parameters, such as $\beta_2$. Furthermore, our numerical continuations show that $\mathbf{SB^{\textit{h}}}$ and $\mathbf{SB^{\ell}}$ are created in a fold bifurcation at $\beta_2 \approx 0.60052$, as the degenerate branching periodic orbit $\mathbf{CSB}$. 

\begin{figure}[t!]
   \centering
   \includegraphics{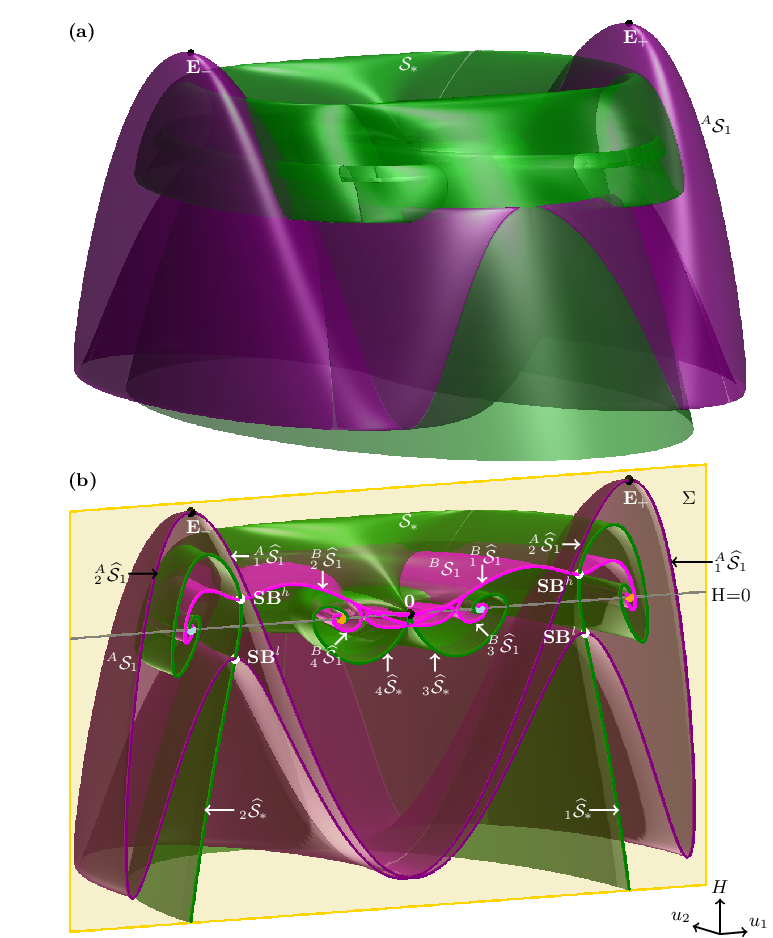}
    \caption{The surfaces $\mathcal{S}_{*}$ (green), ${}^{A}\mathcal{\widehat{S}}_{*}$ (purple) and ${}^{B}\mathcal{\widehat{S}}_{*}$ (magenta) for $\beta_2=0.61$, shown in $(u_1, u_2, H)$-space with $\mathbf{0}$ and $\mathbf{E}_{\pm}$ (black dots). Panel~(a) shows the entire surfaces and panel~(b) is a cutaway view with the section $\Sigma$ defined by $u_2=0$ (beige plane) and the zero-energy level (gray line).}
\label{fig:Sboth_trans}
\end{figure} 

Figure~\ref{fig:Sboth_trans} shows the changed structure of the periodic orbit surfaces past $\mathbf{CSB}$, namely for $\beta_2=0.61$ and again in projection onto $(u_1,u_2,H)$-space with an overall view in panel~(a) and a cutaway view in panel~(b). The surface $\mathcal{S}_{*}$ of $R^*$-symmetric periodic orbits is effectively unchanged, but there are now the new surfaces ${}^{A}\mathcal{\widehat{S}}_{1}$ and ${}^{B}\mathcal{\widehat{S}}_{1}$ of $R_1$-symmetric periodic orbits; compare with Fig.~\ref{fig:Sboth}. The surface ${}^{A}\mathcal{S}_{1}$ contains the equilibria $\mathbf{E}_{\pm}$ and is prominent in panel~(a). However, its actual structure is revealed only in panel~(b), where it is seen to intersect $\mathcal{S}_{*}$ in the symmetry-breaking orbit $\mathbf{SB^{\ell}}$. Note further that ${}^{A}\mathcal{S}_{1}$ intersects the zero-energy level in a single pair of $R_1$-symmetric periodic orbits; in particular, it does not accummulate on the $R_1$-symmetric primary homoclinic orbits and the period of this family of periodic orbits is bounded. The surface ${}^{B}\mathcal{S}_{1}$ is entirely hidden in Fig.~\ref{fig:Sboth_trans}(a), but panel~(b) shows clearly how it intersects $\mathcal{S}_{*}$ in the symmetry-breaking orbit $\mathbf{SB^{h}}$ and accummulates on the pair of the $R_1$-symmetric primary homoclinic orbits. Indeed, the period of these periodic orbits inceases beyond bound in the process. Moroever, the intersection set ${}^{B}\mathcal{\widehat{S}}_{1}$ with the section $\Sigma$ features all the additonal branches that we found for $\mathcal{\widehat{S}}_{1}^{\pm}$ in Sec.~\ref{sec:geomSall}. 

We now consider the consequences of the creation of the symmetry-breaking periodic orbits $\mathbf{SB^{\textit{h}}}$ and $\mathbf{SB^{\ell}}$ at $\beta_2 \approx 0.60052$. Since the degenerate symmetry-breaking periodic orbit $\mathbf{CSB}$ lies in the energy level with $H \approx 0.1652$, both $\mathbf{SB^{\textit{h}}}$ and $\mathbf{SB^{\ell}}$ are initially above the zero-energy level. However, when $\beta_2$ is increased, they separate and at $\beta_2 \approx 0.60064$ the symmetry-breaking periodic orbit $\mathbf{SB^{\ell}}$ crosses the zero-energy level, so that it lies subsequently at negative values of the energy. In fact, Fig.~\ref{fig:Strans}(c) shows the situation for when $\mathbf{SB^{\ell}}$ lies exactly in the zero-energy level, while in Figs.~\ref{fig:Strans}(d) and~\ref{fig:Sboth_trans}(b) for $\beta_2 = 0.61$ it lies well below the zero-energy level. As is illustrated in the enlargement panels~(a2)--(d2), this means that the zero-energy points on the branches ${}_{1}\mathcal{\widehat{S}}_{1}^{-}$ and ${}_{2}\mathcal{\widehat{S}}_{1}^{+}$ disappear when they meet the zero-energy point on ${}_{2}\mathcal{\widehat{S}}_{*}$ at $\beta_2 \approx 0.60064$. Hence, the transition of the symmetry-breaking periodic orbit $\mathbf{SB^{\ell}}$ through the zero-energy level for changing $\beta_2$ induces a generic codimension-one symmetry-breaking bifurcation $\mathbf{SB}$ in the zero-energy level. 

\begin{figure}[t!]
   \centering
   \includegraphics[scale=0.9]{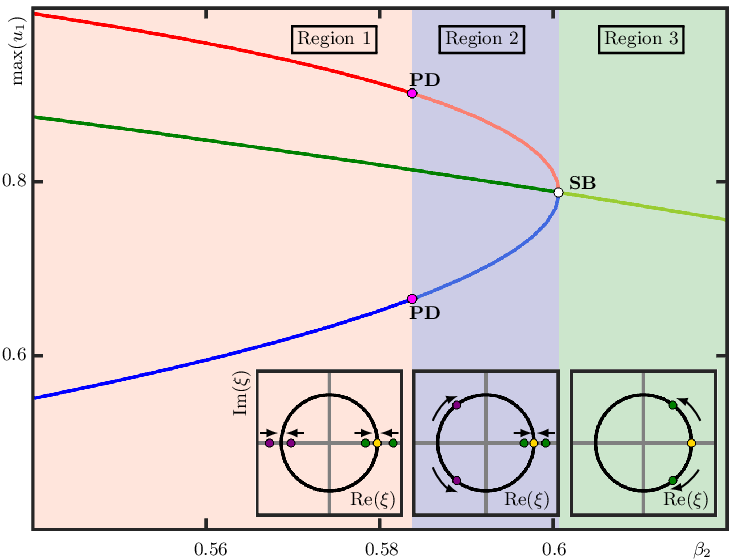}
    \caption{\label{fig:1parbif} One-parameter bifurcation diagram in $\beta_2$ of the periodic orbits in the zero-energy level near the symmetry-breaking bifurcation $\mathbf{SB}$, where the branches of $R^{*}$-symmetric and $R_1$-symmetric periodic orbits are represented by $\max{u_1}$ and are coloured as in Fig.~\ref{fig:Strans}(a2)-(d2). The insets show in the coloured regions 1 to 3 show the unit circle in the complex plane with the trivial multiplier $1$ (yellow dot), and the two non-trivial Floquet multipliers of the $R^{*}$-symmetric periodic orbit (green dots) and of the pair of $R_1$-symmetric periodic orbits (purple); the the arrows indicate how the respective Floquet multipliers change for increasing $\beta_2$.} 
\end{figure} 

Figure~\ref{fig:1parbif} shows the one-parameter bifurcation diagram in $\beta_2$ with the symmetry-breaking bifurcation $\mathbf{SB}$ at $\beta_2 \approx 0.60064$, where a pair branches of $R_1$-symmetric periodic orbits bifurcates from the branch of $R^{*}$-symmetric periodic orbits, which exists for $\beta_2 < 0.60064$. Notice that the branches of $R_1$-symmetric periodic orbits meet at the point $\mathbf{SB}$ to form a single smooth curve. At the point $\mathbf{SB}$ the point $\mathbf{SB}$ there is a qualitative change of the Floquet multipliers of the $R^{*}$-symmetric periodic orbits. In divergence-free vector fields, such as~\eqref{eq:4Dode}, the product of the Floquet multipliers is equal to one \cite{Kent_1992}, which means that there are always two trivial Floquet multipliers equal to $1$ as well as two non-trivial multipliers whose product is $1$. 
The insets of Fig.~\ref{fig:1parbif} show the Floquet multipliers of all existing periodic orbits in the highlighted regions 1--3. At the point $\mathbf{SB}$ the two non-trivial Floquet multipliers of the $R^{*}$-symmetric periodic orbits in the zero-energy level coincide at $1$, which characterises a symmetry-breaking bifurcation of a Hamiltonian system \cite{Kent_1992}. At this point, these periodic orbits transition from elliptic (with complex conjugate non-trivial eigenvalues on the unit circle) in regions 1 and 2 to hyperbolic (with real non-trivial eigenvalues) in region 3. 

In region 2 the $R_1$-symmetric periodic orbits in the zero-energy level are elliptic; as $\beta_2$ is increased, the complex conjugate pair of their Floquet multipliers moves along the unit circle and meets at $1$ when $\mathbf{SB}$ at $\beta_2 \approx 0.60064$ is reached. However, as $\beta_2$ is decreased in region 2, these Floquet multipliers meet at $-1$ when $\beta_2 \approx 0.58373$. This point, labeled $\mathbf{PD}$ in Fig.~\ref{fig:1parbif}, is a period-doubling bifurcation \cite{Kent_1992} of the Hamiltonian system~\eqref{eq:4Dode}. As is shown by the insets, as a result the $R_1$-symmetric periodic orbits are hyperbolic in region 1. We remark that, since their non-trivial real Floquet multipliers are negative, the $R_1$-symmetric periodic orbits are non-orientable in region 1, that is, have stable and unstable manifolds that are M{\"o}bius bands locally near these periodic orbits. The $R^{*}$-symmetric periodic orbits in the zero-energy, on the other hand, have positive non-trivial real Floquet multipliers in regions 1 and 2; hence, they are orientable, that is, their stable and unstable manifolds are locally topological cylinders \cite{Kent_1992}.

\section{Period-$k$ multiplying bifurcations of $R^{*}$-symmetric periodic orbits}
\label{sec:periodmult}

The changes to the structure of periodic orbits in the last section concerned chiefly the $R_1$-symmetric periodic orbits, while the surface $\mathcal{S}_{*}$ of the $R^*$-symmetric periodic orbits remained unchanged. As we show now, when $\beta_2$ is increasing further, $\mathcal{S}_{*}$ interacts with other surfaces of $R^{*}$-symmetric periodic orbits in a series of bifurcations known as period-$k$ multiplying bifurcations \cite{GREENE1981468,mackay1982renormalisation,Kent_1992,webster_2003}. These can be understood by considering the first return map to a local Poincar{\'e} section, such as the section $\Sigma$ defined by $u_2=0$. Since system~\eqref{eq:4Dode} is Hamiltonian, this map is volume preserving and, moreover, its restriction to a particular energy level is a two-dimensional area-preserving diffeomorphism. MacKay \emph{et al.} \cite{GREENE1981468,mackay1982renormalisation} studied period-$k$ multiplying bifurcations of fixed points in this class of maps for $k \geq 3$. Such a bifurcation occurs when an elliptic fixed point has a pair of complex conjugate eigenvalues with rational rotation number, which are, hence, of the form $e^{\pm 2\pi i\frac{p}{k}}$ for some $p\in \mathbb{N}$ that is incommensurate with $k$.

The unfolding for $k=3$ and $k=4$ differs from that for $k \geq 5$, and these two cases from \cite{mackay1982renormalisation,Kent_1992} are shown schematically in Fig.~\ref{fig:kmult_sketch} as adapted to the setting here. Panel~(a) illustrates the case $k=3,4$, where a branch of elliptic periodic orbits and a branch of hyperbolic periodic orbits of $k$ times the period meet and cross; since this is reminiscent of a transcritical bifurcation, we refer to this period-$k$ multiplying bifurcation and the corresponding periodic orbit as $\mathbf{T_k}$. Moreover, the hyperbolic periodic orbits become elliptic at a nearby saddle-node bifurcation $\mathbf{SN}$. Figure~\ref{fig:kmult_sketch}(b) shows the case $k \geq 5$, where two simultaneously existing branch of periodic orbits of $k$ times the period, one hyperbolic and the other elliptic,emerge from (or disappear at) the period-$k$ multiplying bifurcation denoted $\mathbf{B_k}$. We remark that the cases $k=1$ and $k=2$ also exist, but are not referred to as period-$k$ multiplying bifurcations: for $k=1$ there is a double eigenvalue $+1$, which is the symmetry-breaking bifurcation, and for $k=2$ there is a double eigenvalue $-1$, which is the period-doubling bifurcation. 

Period-$k$ multiplying bifurcations have been observed in a number of divergence-free systems; for example, in the Michelson system, which is of third-order, reversible and volume-preserving \cite{Kent_1992,webster_2003}. In system~\eqref{eq:4Dode} the eigenvalue condition concerns the non-trivial Floquet multipliers of the respective periodic orbits. Indeed, we already identified in Fig.~\ref{fig:1parbif} the cases of symmetry-breaking and period-doubling bifurcations. Moreover, we saw that the Floquet multipliers of elliptic $R^*$-symmetric periodic as well as $R_1$-symmetric periodic orbits in the zero-energy surface move over the unit circle when $\beta_2$ is changed and will, hence, undergo period-$k$ multiplying bifurcations. As we will show, these bifurcations are induced by the creation in pairs of period-$k$ multiplying periodic orbits, which exist for fixed $\beta_2$ at certain levels of the Hamiltonian energy $H$. This scenario is indeed completely analogous to what we found for the symmetry-breaking bifurcation in Sec.~\ref{sec:symmbreak}. 

More specifically, the elliptic single-loop $R^*$-symmetric periodic orbits on the surface $\mathcal{S}_{*}$ we consider in this section have period-$k$ multiplying bifurcations for any odd $k$. We first discuss the transcritical-type period-$3$ multiplying bifurcation and its consequences, and then subsequent period-$k$ multiplying bifurcations for $k = 5,7,...$ .

\begin{figure}[t!]
   \centering
   \includegraphics{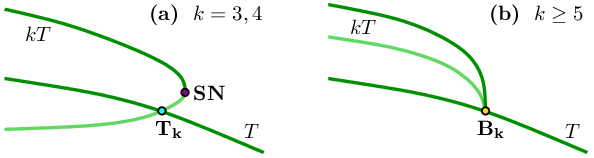}
    \caption{Schematic one-parameter bifurcation diagrams of the period-$k$ multiplying bifurcation of elliptic periodic orbits, adapted from \cite{mackay1982renormalisation,Kent_1992}. Panel~(a) shows the transcritical-type case for $k=3,4$ with the points $\mathbf{T_k}$ and $\mathbf{SN}$, and panel~(b) the general case for $k\geq5$ with the bifurcation referred to as $\mathbf{B_k}$. The horizontal represents the Hamiltonian energy $H$, branches of elliptic periodic orbits are dark green, those of hyperbolic periodic orbits light green, and their period $T$ or $kT$ is indicated.} 
       \label{fig:kmult_sketch}
\end{figure} 

\subsection{Transcritical-type period-$3$ multiplying bifurcation}
\label{sec:per3mult}

\begin{figure}[t!]
   \centering
   \includegraphics[width=12cm]{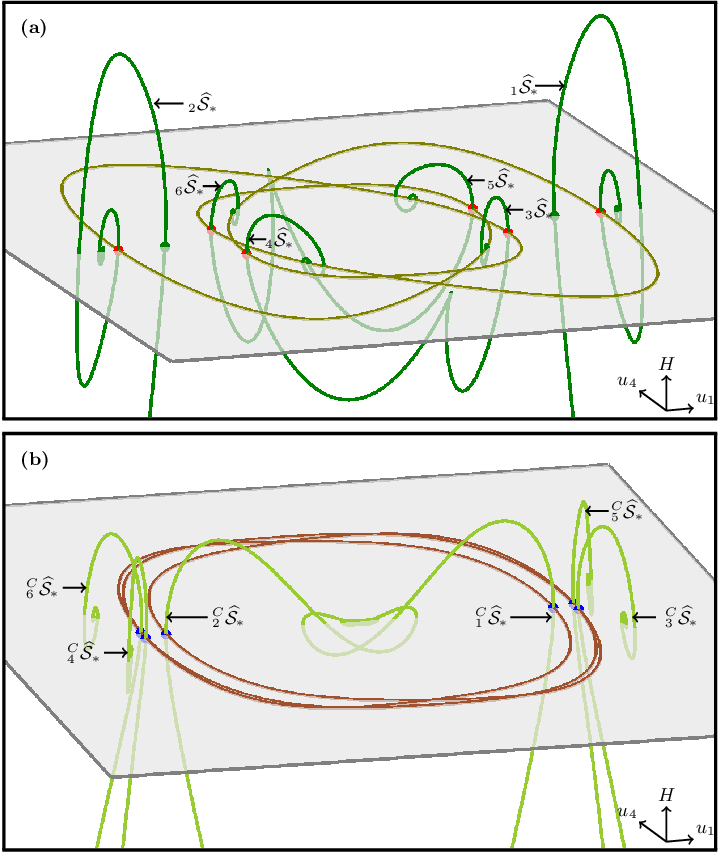}
       \caption{\label{fig:per3int} Intersection sets with the section $\Sigma$ for  $\beta_2=0.63$ of the surfaces $\mathcal{S}_{*}$ (dark green) in panel~(a) and ${}^{C}\mathcal{S}_{*}$ (light green) in panel~(b), shown in projection onto $(u_1,u_4H)$-space. Also shown is the zero-energy level (gray plane) and the respective period-one and period-three periodic orbit in it, with points in $\Sigma$ marked by green, red and blue dots.} 
\end{figure} 

Numerical continuation shows that at $\beta_2=0.6397$ the single-loop $R^{*}$-symmetric periodic orbits on $\mathcal{S}_{*}$ features a degenerate transcritical-type period-$3$ multiplying periodic orbit $\mathbf{CT_3}$ in the energy level with $H\approx 0.1245$. This bifurcation marks the creation of a pair $\mathbf{T_3^{\textit{h}}}$ and $\mathbf{T_3^{\ell}}$ of period-$3$ multiplying periodic orbits with higher and lower $H$, respectively. The period-three orbits that connect at $\mathbf{CT_3}$ form a surface we denote ${}^{C}\mathcal{S}_{*}$. Understanding and illustrating this and subsequent period-$k$ multiplying bifurcations requires a three-dimensional view of the section $\Sigma$ (defined by $u_2 = 0$); a two-dimensional projection as in Fig.~\ref{fig:Strans} is no longer sufficient, because some curves in the respective intersections sets have the same projections onto the $(u_1,H)$-plane. To this end, we represent $\Sigma$ by $(u_1,u_4,H)$-space from now on.

Figure~\ref{fig:per3int} illustrates in this way the intersection sets of the surfaces $\mathcal{S}_{*}$ and ${}^{C}\mathcal{S}_{*}$ for $\beta_2=0.63$ before $\mathbf{CT_3}$. Also shown here is the zero-energy level (now a plane) and respective $R^{*}$-symmetric periodic orbits in it. Panel~(a) shows the $H$-parametrised intersections set of $\mathcal{S}_{*}$ with the six labeled branches $_{i}\mathcal{\widehat{S}_{*}}$ with $i=1,...,6$ introduced in Fig.~\ref{fig:Sboth}. Note in Fig.~\ref{fig:per3int}(a) that the first additonal pairs of branches $_{3}\mathcal{\widehat{S}_{*}}$, $_{5}\mathcal{\widehat{S}_{*}}$ and $_{4}\mathcal{\widehat{S}_{*}}$, $_{6}\mathcal{\widehat{S}_{*}}$ are now distinguished in $(u_1,u_4,H)$-space. The first points of the branches $_{1}\mathcal{\widehat{S}_{*}}$ and $_{2}\mathcal{\widehat{S}_{*}}$ in the zero-energy level are the intersection points with $\Sigma$ of the shown single-loop periodic orbit, which is the continuation of $\Gamma_{*}$. The second point in the zero-energy level on each of these two branches, together with the first points of the four additional branches, are the intersection points with $\Sigma$ of the triple-loop periodic orbit that is also shown. Fig.~\ref{fig:per3int}(b) shows in the same way six labeled branches ${}^{C}_{i}\mathcal{\widehat{S}_{*}}$ with $i=1,...,6$ of the intersection set of ${}^{C}\mathcal{S}_{*}$. Near their first intersection points with the zero-energy level, they correspond to triple-loop periodic orbits as the one shown with $H=0$.

\begin{figure}[t!]
   \centering
   \includegraphics[width=12cm]{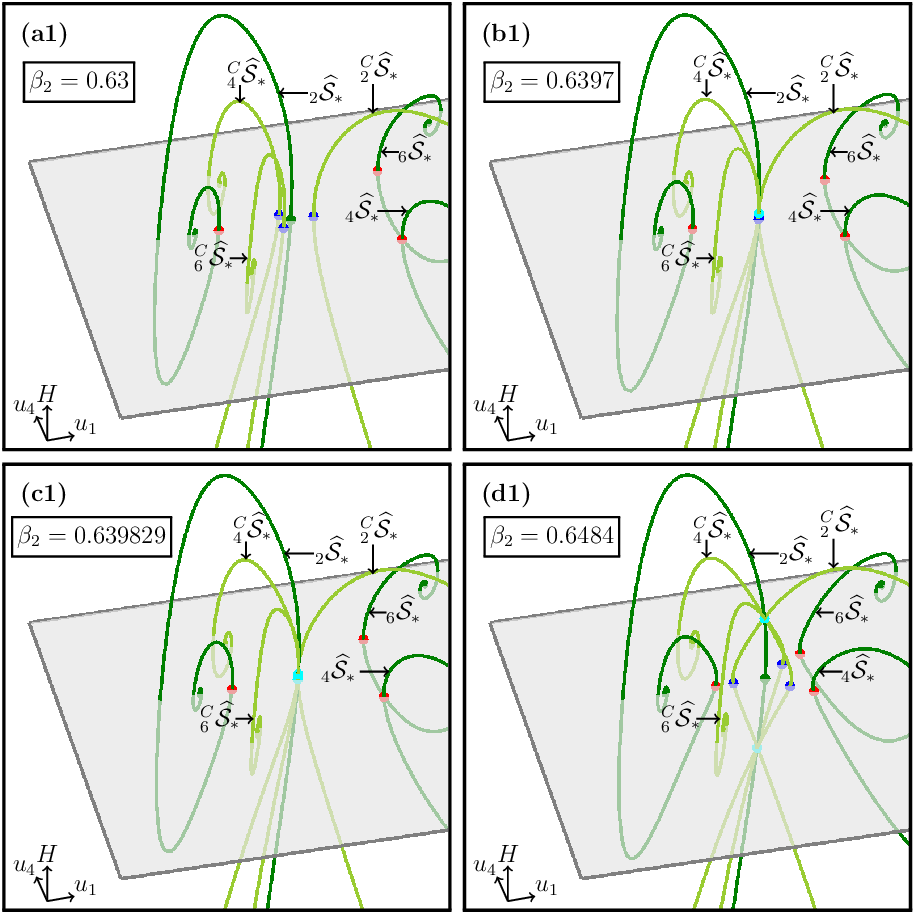}
       \caption{\label{fig:per3trans} Transition through and past $\mathbf{CT_3}$, illustrated in $\Sigma$ near $_{2}\mathcal{\widehat{S}_{*}}$ (dark green) in projection onto $(u_1,u_4H)$-space. Panels~(a1)--(f1) give an overall view of insection sets of the surfaces ${}^{C}\mathcal{S}_{*}$ (light green), ${}^{D}\mathcal{S}_{*}$ (steel blue), and ${}^{E}\mathcal{S}_{*}$ (mid-green) for the stated $\beta_2$-values, and panels~(a2)--(f2) are enlargements near the zero-energy level. Also shown are the respective intersection points (green, red and blue dots) in the zero-energy level (grey plane), as well as the points $\mathbf{CT_3}$, $\mathbf{T_{3}^{\textit{h}}}$ and $\mathbf{T_{3}^{\ell}}$ (cyan dots), and $\mathbf{CSN}$, $\mathbf{SN^{\textit{h}}}$ and  $\mathbf{SN^{\ell}}$ (purple dots).} 
\end{figure} 
\addtocounter{figure}{-1}
\begin{figure}[t!]
   \centering
      \includegraphics[width=12cm]{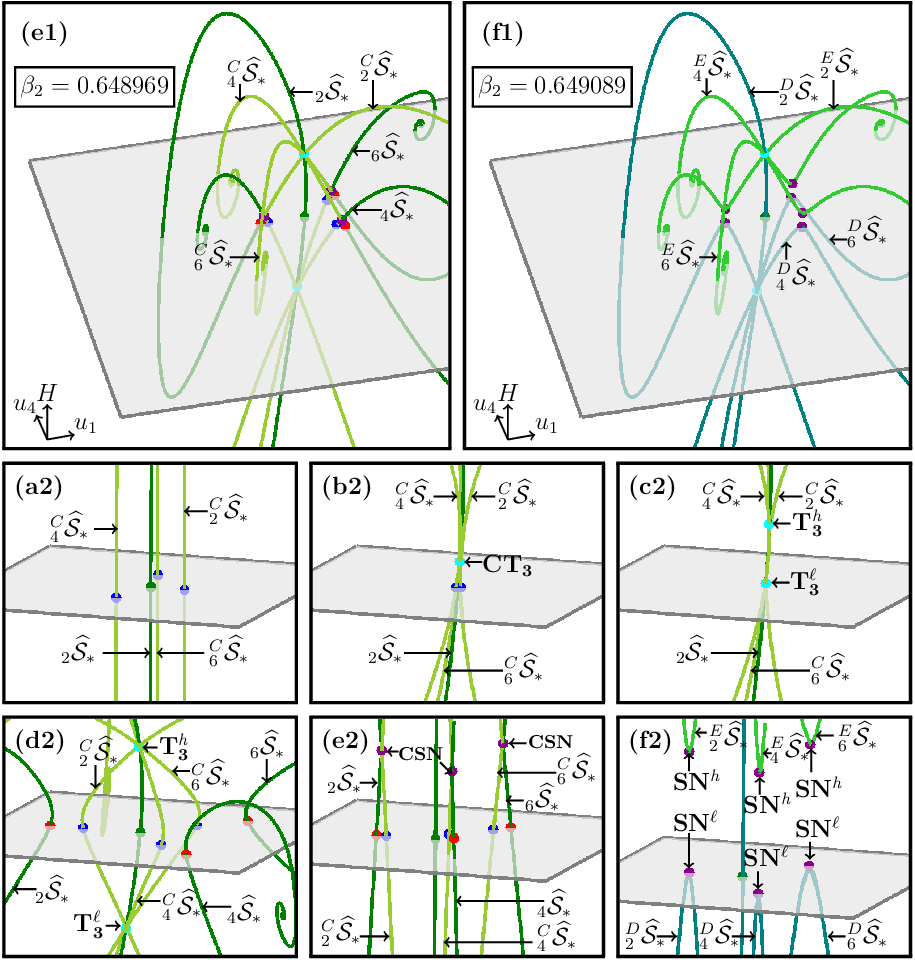}
       \caption{Continued.} 
\end{figure} 

The different branches in Fig.~\ref{fig:per3int} are well separated, but they approach each other and then interact as $\beta_2$ is increased. This is illustrated in Fig.~\ref{fig:per3trans}, where we show changes to the local arrangement of the branches $_{i}\mathcal{\widehat{S}_{*}}$ and ${}^{C}_{i}\mathcal{\widehat{S}_{*}}$ with even $i$. Panel~(a1) is again for $\beta_2=0.63$ and shows the overall situation. As the enlargement panel~(a2) clearly shows, the branches ${}^{C}_{2}\mathcal{\widehat{S}}_{*}$, ${}^{C}_{4}\mathcal{\widehat{S}}_{*}$ and ${}^{C}_{6}\mathcal{\widehat{S}}_{*}$ ``surround'' the branch $_{2}\mathcal{\widehat{S}}_{*}$ near the single-loop periodic orbit in the zero-energy level. These branches are, in turn, surrounded in panel~(a1) by ${}_{2}\mathcal{\widehat{S}}_{*}$, $_{4}\mathcal{\widehat{S}}_{*}$ and $\mathcal{\widehat{S}}_{*}$ near their triple-loop periodic orbit in the zero-energy level; however, these branches are too far away to be visible in panel~(a2).

Figure~\ref{fig:per3trans}(b1) for $\beta_2=0.6397$ shows the approximate moment when the branches ${}^{C}_{2}\mathcal{\widehat{S}}_{*}$, ${}^{C}_{4}\mathcal{\widehat{S}}_{*}$ and ${}^{C}_{6}\mathcal{\widehat{S}}_{*}$ meet $_{2}\mathcal{\widehat{S}}_{*}$ tangentially at the $\mathbf{CT_3}$, which  lies above the zero-energy level at $H\approx0.12451$; see panel~(b2). Hence, the periodic orbits of ${}^{C}\mathcal{S}_{*}$ become the triple cover of the single-loop periodic of $\mathcal{S}_{*}$ in this energy level. The point $\mathbf{CT_3}$ creates the two transcritical-type period-$3$ multiplying periodic orbits $\mathbf{T_{3}^{\textit{h}}}$ and $\mathbf{T_{3}^{\ell}}$ for increasing $\beta_2$, which are both initially above the zero-energy level. At $\beta_2\approx0.639829$, as in Fig.~\ref{fig:per3trans}(c2), the lower point $\mathbf{T_{3}^{\ell}}$ crosses the zero-energy level into the region of negative $H$. Hence, at this value of $\beta_2$ there is an induced transcritical-type period 3-multiplying bifurcation in the zero-energy level. For even larger $\beta_2$, as in Fig.~\ref{fig:per3trans}(d1), the point $\mathbf{T_{3}^{\ell}}$ is well below the zero-energy level, while $\mathbf{T_{3}^{\textit{h}}}$ has moved further up in $H$. Notice here that the branches 
${}^{C}_{2}\mathcal{\widehat{S}}_{*}$, ${}^{C}_{4}\mathcal{\widehat{S}}_{*}$ and ${}^{C}_{6}\mathcal{\widehat{S}}_{*}$ each intersect $_{2}\mathcal{\widehat{S}}_{*}$ transversally at the points $\mathbf{T_{3}^{\textit{h}}}$ and $\mathbf{T_{3}^{\ell}}$. Moreover their intersection points with the zero-energy level are now further away from that of the branch $_{2}\mathcal{\widehat{S}}_{*}$, while the surrounding branches ${}^{C}_{2}\mathcal{\widehat{S}}_{*}$, ${}^{C}_{4}\mathcal{\widehat{S}}_{*}$ and ${}^{C}_{6}\mathcal{\widehat{S}}_{*}$ are now closer; see also panel~(d2). 

As Fig.~\ref{fig:per3trans}(e1) and the enlargement panel~(e2) show, at $\beta_2 \approx 0.648969$ the intersection branches ${}^{C}_{2}\mathcal{\widehat{S}}_{*}$, ${}^{C}_{4}\mathcal{\widehat{S}}_{*}$ and ${}^{C}_{6}\mathcal{\widehat{S}}_{*}$ meet 
$_{2}\mathcal{\widehat{S}}_{*}$, $_{4}\mathcal{\widehat{S}}_{*}$ and $_{6}\mathcal{\widehat{S}}_{*}$, respectively, tangentially at three points labelled $\mathbf{CSN}$ in the energy level with $H\approx0.12451$. This results in different parts of these branches connecting differently for larger $\beta_2$, which gives rise to the new intersection sets ${}^{D}\mathcal{\widehat{S}}_{*}$ and ${}^{E}\mathcal{\widehat{S}}_{*}$ shown in panels~(f1) and~(f2) of Fig.~\ref{fig:per3trans}. More specifically, the special periodic orbit $\mathbf{CSN}$ creates two saddle-node periodic orbits, namely a local maximum $\mathbf{SN^{\textit{h}}}$ and a local minimum $\mathbf{SN^{\ell}}$ of $H$, respectively, which are initially above the zero-energy level. Figure~\ref{fig:per3trans}(f2) actually shows the moment that $\mathbf{SN^{\ell}}$ crosses the zero-energy surface at $\beta_2\approx0.649089$; hence, this is an induced saddle-node bifurcation of periodic orbits in the zero energy surface. 

In fact, past the point $\mathbf{CSN}$, as in Fig.~\ref{fig:per3trans}(f2), we find exactly two copies of the situation sketched in Fig.~\ref{fig:kmult_sketch}(a), with the points $\mathbf{T_{3}^{\textit{h}}}$ and $\mathbf{SN^{\textit{h}}}$, and $\mathbf{T_{3}^{\ell}}$ and $\mathbf{SN^{\ell}}$, respectively. The difference is that the theory considers the third iterate of the return map and, hence, the sketch show a single and not three branches. Indeed, the nearby degenerate saddle-node bifurcation $\mathbf{CSN}$, rather than initial bifurcation $\mathbf{CT_{3}}$ is responsible for the vastly changed geometry of the periodic orbit structure in Fig.~\ref{fig:per3trans}(f2). Notice, in particular, that the surface ${}^{D}\mathcal{S}_{*}$ does not accumulate on homoclinic orbits, while the surface ${}^{E}\mathcal{S}_{*}$ does.

\subsection{Nested period-k multiplying bifurcations for $k\geq5$}
\label{sec:nested}

\begin{figure}[t!]
   \centering
   \includegraphics[width=12cm]{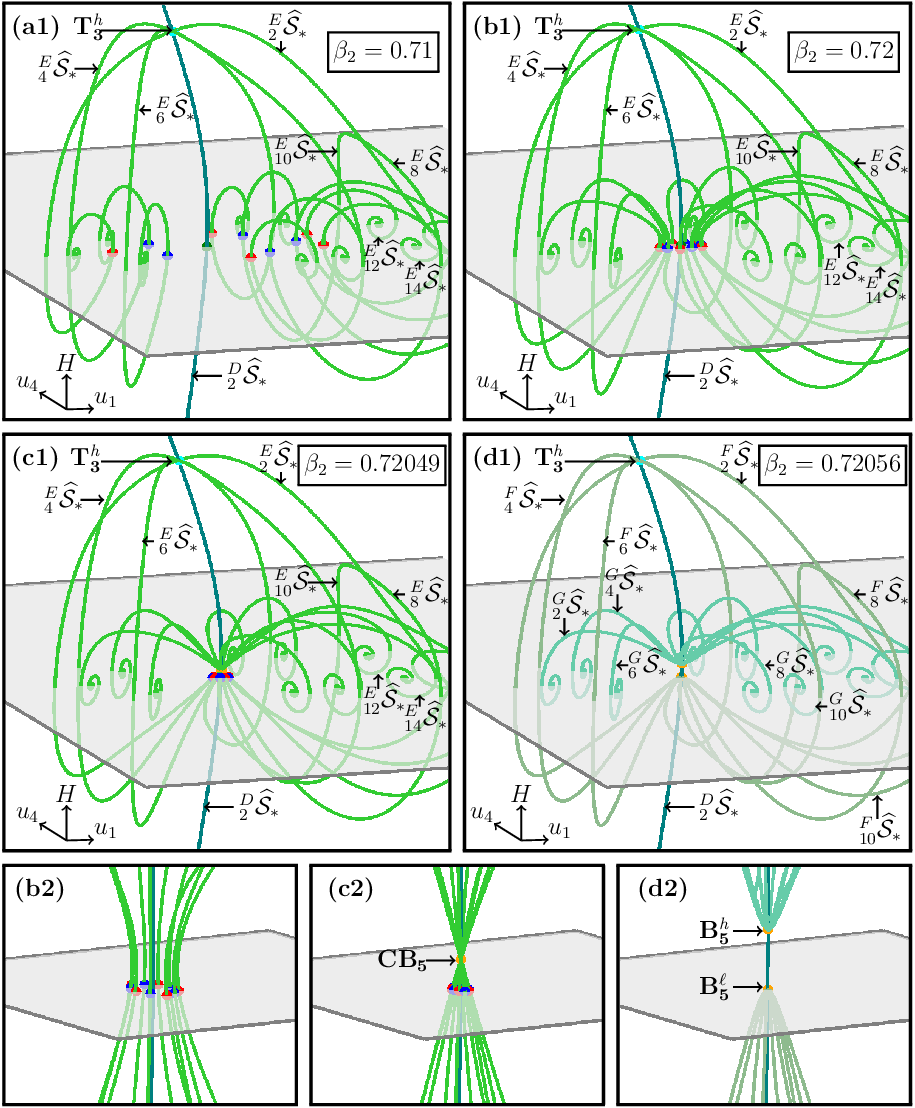}
       \caption{\label{fig:per5trans} Transition through and past $\mathbf{CB_5}$, illustrated in $\Sigma$ near ${}^{D}_{2}\mathcal{\widehat{S}_{*}}$ (steel blue) in projection onto $(u_1,u_4H)$-space. Panels~(a1)--(d1) give an overall view of the involved (again even) branches of the insection set of the surface ${}^{E}\mathcal{S}_{*}$ (mid green) for the stated $\beta_2$-values, and panels~(b2)--(d2) are enlargements near the zero-energy level. Also shown are the respective intersection points (green, red and blue dots) in the zero-energy level (grey plane), as well as the points $\mathbf{CB_5}$, $\mathbf{B_{5}^{\textit{h}}}$ and $\mathbf{B_{5}^{\ell}}$ (orange dots) and $\mathbf{T_{3}^{\textit{h}}}$ (cyan dot).}  
\end{figure} 

Figure~\ref{fig:per5trans} illustrates that the same central branch, now labeled ${}^{D}_{2}\mathcal{\widehat{S}}_{*}$ because it is part of the new surface ${}^{D}\mathcal{S}_{*}$, has a degenerate period-$5$ multiplying bifurcation $\mathbf{CB_5}$ at $\beta_2 \approx 0.720485$, which creates two points $\mathbf{B_5^{\textit{h}}}$ and $\mathbf{B_5^{\ell}}$, with the latter moving through  zero-energy level for slightly higher $\beta_2$. For orientation, we show here also the transcritical-type period-$3$ multiplying bifurcation $\mathbf{T_{3}^{\textit{h}}}$ above the zero-energy level. In Fig.~\ref{fig:per5trans}(a) before this transition, two sets of five intersection points in the zero-energy surface on the respective shown branches of the intersection set ${}^{E}_{2}\mathcal{\widehat{S}}_{*}$ correspond to a pair of five-loop $R^{*}$-symmetric periodic orbits. These points and the corresponding parts of the branches of ${}^{E}_{2}\mathcal{\widehat{S}}_{*}$ move closer to the branch ${}^{D}_{2}\mathcal{\widehat{S}}_{*}$ see Fig.~\ref{fig:per5trans}(b1) and the enlargement panel~(b2). At $\beta_2 \approx 0.720485$ these branches meet tangentially at the point $\mathbf{CB_5}$, which lies in the energy level with $H\approx0.057063$, as is illustrated in panels~(c1) and (c2). The point $\mathbf{CB_5}$ represents the degenerate period-$5$ multiplying periodic orbit, which creates, for increasing $\beta_2$, two period-$5$ multiplying periodic orbits $\mathbf{B_5^{\textit{h}}}$ and $\mathbf{B_5^{\ell}}$, of which $\mathbf{B_5^{\textit{h}}}$ and $\mathbf{B_5^{\ell}}$ moves through $H=0$ at $\beta_2=0.720559$. This is, hence, the moment of when there is an induced period-$5$ multiplying bifurcation in the the zero-energy level, which is the situation shown in Fig.~\ref{fig:per5trans}(d1) and~(d2). Note that the local picture near the points $\mathbf{B_5^{\textit{h}}}$ and $\mathbf{B_5^{\ell}}$ is exactly as sketched in Fig.~\ref{fig:kmult_sketch}(b) for the fifth iterate of the return map.

The associated local re-arrangement of the $R^{*}$-symmetric periodic orbits on ${}^{E}_{2}\mathcal{S}_{*}$ past the point $\mathbf{CB_5}$ leads to the creation of the new surfaces ${}^{F}\mathcal{S}_{*}$ and ${}^{G}\mathcal{S}_{*}$. The branches of the intersection set ${}^{F}\mathcal{\widehat{S}}_{*}$ in Fig.~\ref{fig:per5trans}(d1) have the period-$3$ multiplying periodic orbit $\mathbf{T_{3}^{\textit{h}}}$ on ${}^{D}_{2}\mathcal{\widehat{S}}_{*}$ as their maxima, and the period-$5$ multiplying periodic orbit $\mathbf{B_{5}^{\ell}}$ also on ${}^{D}_{2}\mathcal{\widehat{S}}_{*}$ as their minima. In other words, the surface ${}^{F}_{2}\mathcal{S}_{*}$ consists of periodic orbits with bounded period. The surface ${}^{G}\mathcal{S}_{*}$, on the other hand, accumulates on homoclinic orbits in the zero-enery surfaces, that is, its intersection set ${}^{F}\mathcal{\widehat{S}}_{*}$ has spiraling branches and the period is unbounded.

\begin{figure}[t!]
   \centering
   \includegraphics[width=12cm]{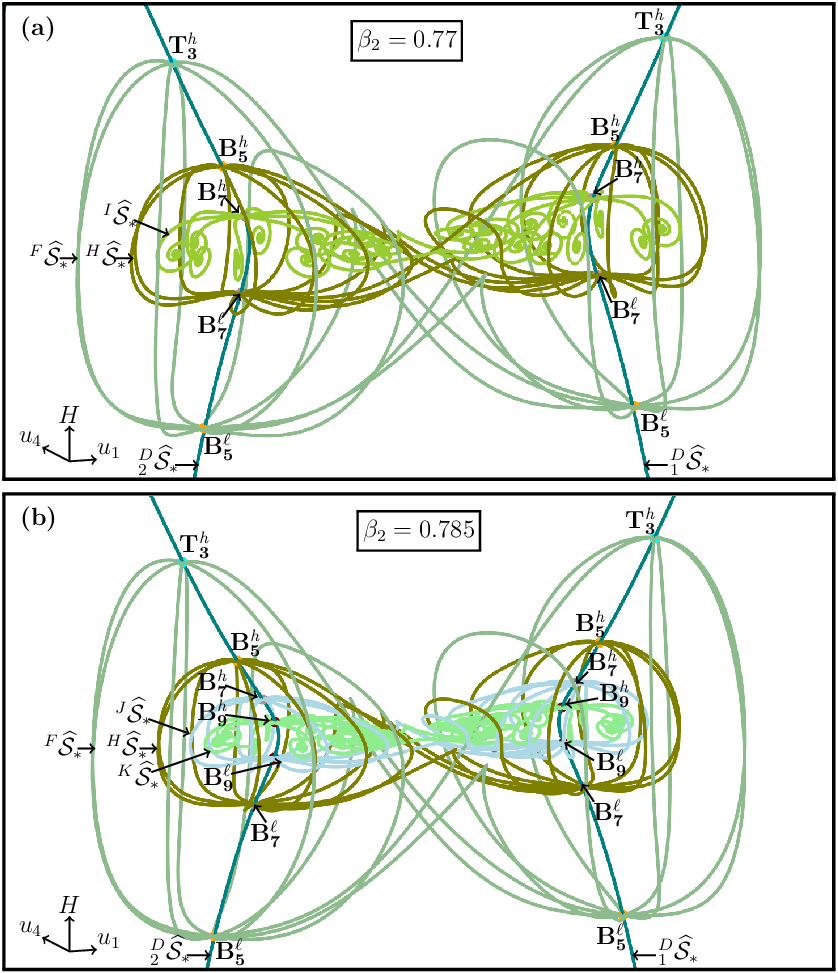}
       \caption{\label{fig:global7and9} Global view in $(u_1,u_4,H)$-space of how the relevant intersections sets with $\Sigma$ of surfaces of $R^{*}$-symmetric periodic orbits connect at period-$k$ multiplying bifurcations on ${}^{D}\mathcal{\widehat{S}}_{*}$ for $\beta_2=0.77$ after $\mathbf{CB_{7}}$ in panel~(a), and for $\beta_2=0.785$ after $\mathbf{CB_{9}}$ in panel~(b).}
\end{figure} 

Figure~\ref{fig:global7and9} shows that there is a sequence of further period-$k$ multiplying bifurcations for odd $k \geq 7$ that each lead to similar bfurcations of the periodic orbit structure. Here we show a global view in $(u_1,u_4,H)$-space with the intersection sets of all the respective surfaces and period-$k$ multiplying periodic orbits. Specifically, panel~(a) shows the situation after the bifurcation $\mathbf{CB_{7}}$, which created the pair of period-$7$ multiplying periodic orbits $\mathbf{B_{7}^{\textit{h}}}$ and $\mathbf{B_{7}^{\ell}}$. As a result, the intersection set ${}^{G}\mathcal{\widehat{S}}_{*}$ has bifurcated into the set ${}^{H}\mathcal{\widehat{S}}_{*}$, which connects $\mathbf{B_{5}^{\textit{h}}}$ and $\mathbf{B_{7}^{\ell}}$ and has bounded period, and ${}^{I}\mathcal{\widehat{S}}_{*}$ whose branches spiral into the intersection points of the $R_1$-symmetric and $R_2$-symmetric primary homoclinic orbits. Figure~\ref{fig:global7and9}(b) shows that this process repeats when $\mathbf{CB_{9}}$ is passed: now ${}^{H}\mathcal{\widehat{S}}_{*}$ bifurcated into the new intersection sets ${}^{J}\mathcal{\widehat{S}}_{*}$, which connects $\mathbf{B_{7}^{\textit{h}}}$ and $\mathbf{B_{9}^{\ell}}$, and ${}^{K}\mathcal{\widehat{S}}_{*}$ with spiraling branches.

Overall, Figure~\ref{fig:global7and9} shows how the sequence of degenerate period-$k$ multiplying bifurcations $\mathbf{CB_{k}}$ with odd $k$ creates nested subsurfaces of $R^{*}$-symmetric periodic orbits with bounded period --- each of which connecting the period-$k$ multiplying bifurcations $\mathbf{B_{k-2}^{\textit{h}}}$ and $\mathbf{B_{k}^{\ell}}$. In the limit $k \to \infty$, this sequence itself accumulates on the $R_1$-symmetric and $R_2$-symmetric primary homoclinic orbits in the zero-energy level. Our numerical evidence shows that the special periodic orbits  $\mathbf{CB_k}$ for $k=3,5,7$ all have $H>0$ (for the chosen values of the other parameters); in fact, the sequence $H(\mathbf{CB_{k}})$ is strictly decreasing and converges to $0$ in the limit. As a consequence, the corresponding points $\mathbf{B_{k}^{\ell}}$ are crossing the zero-energy level for $\beta_2$-values ever closer to those where one finds $\mathbf{CB_{k}}$.

\section{Symmetry-breaking period-$k$ multiplying bifurcations}
\label{sec:nonsymm}

\begin{figure}[t!]
   \centering
   \includegraphics[width=12cm]{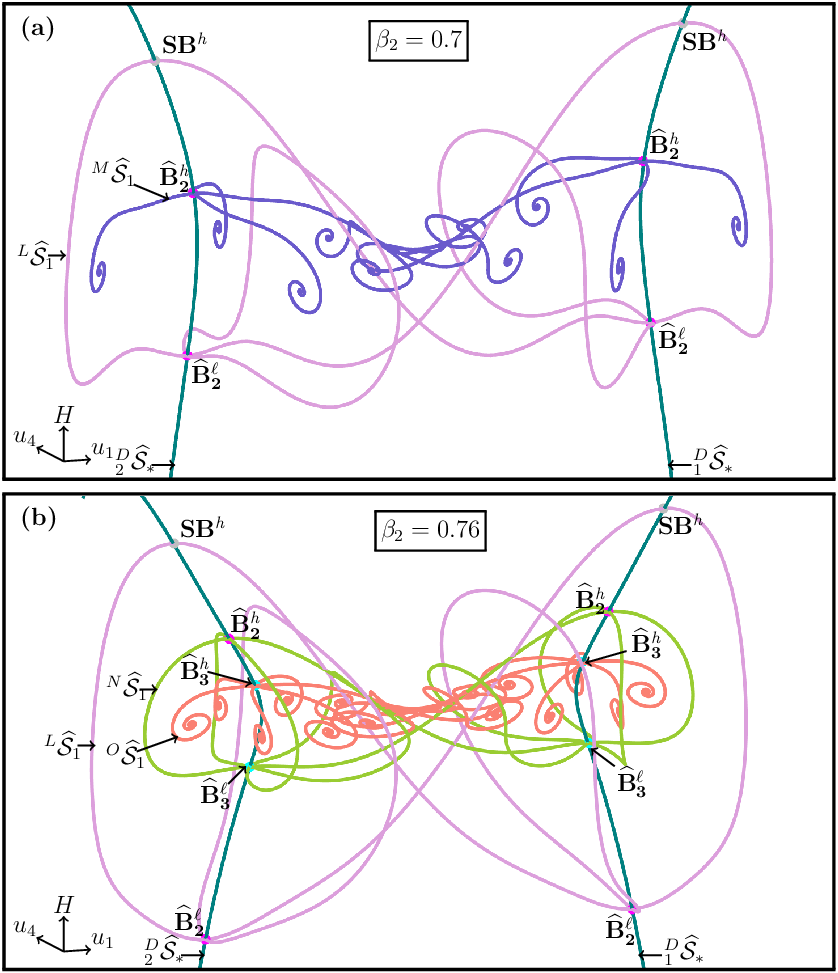}
       \caption{\label{fig:hatglobal2and3} Global view in $(u_1,u_4,H)$-space of how the intersection set ${}^{D}\mathcal{\widehat{S}}_{*}$ of $R^{*}$-symmetric periodic orbits connects with the intersections sets of $R_1$-symmetric periodic orbits at symmetry-breaking period-$k$ multiplying bifurcations for $\beta_2=0.7$ after $\mathbf{\widehat{CB}_{2}}$ in panel~(a), and for $\beta_2=0.76$ after $\mathbf{\widehat{CB}_{3}}$ in panel~(b).} 
\end{figure} 

The period-$k$ multiplying bifurcations $\mathbf{B_{k}}$ of $R^{*}$-symmetric periodic orbits we studied so far concerned bifurcating branches of period-$k$ periodic orbits that are also $R^{*}$-symmetric. We now show that there is also a different and new type, which we refer to as symmetry-breaking period-$k$ multiplying bifurcations of $R^{*}$-symmetric periodic orbits. Their distinct feature is that from them bifurcate branches/surfaces of only $R_{1}$-symmetric period-$k$ periodic orbits, as well as of only $R_{2}$-symmetric period-$k$ periodic orbits; indeed these come in pairs as each others image under the respective other reversing symmetry. More specifically, we identify a sequence of special periodic orbits $\mathbf{\widehat{CB}_{k}}$ at specific and increasing values of $\beta_2$, which creates pairs of symmetry-breaking period-$k$ multiplying bifurcations $\mathbf{\widehat{B}_{k}^{\textit{h}}}$ and $\mathbf{\widehat{B}_{k}^{\ell}}$ for any $k \geq 2$ (not just odd $k$). Analoguous to what we found in the last section, the points $\mathbf{\widehat{CB}_{k}}$ all have values of $H$ that are positive and converge monotonically to zero as $k \to \infty$; again the points $\mathbf{\widehat{B}_{k}^{\ell}}$ then cross the zero-energy surface and, thus, generate symmetry-breaking period-$k$ multiplying bifurcations of periodic orbits with $H=0$.

The symmetry-breaking period-$k$ multiplying bifurcations also change the geometry of surfaces of periodic orbits, with the difference that the periodic orbits that are involved have less symmetry. We do not go into quite as much detail as we did in Sec.~\ref{sec:periodmult} and only show how the surfaces are organised after the first two degenerate symmetry-breaking period-$k$ multiplying bifurcations $\mathbf{\widehat{CB}_{2}}$ and $\mathbf{\widehat{CB}_{3}}$. More specifically, we only consider here the geometrical changes that happen to the surfaces of $R_1$-symmetric periodic orbits. Figure~\ref{fig:hatglobal2and3} shows the organisation of the relevant surfaces of periodic orbits on the level of their intersections sets with the section $\Sigma$, represented in the style of Figure~\ref{fig:global7and9} by the projection onto $(u_1,u_4,H)$-space.

Figure~\ref{fig:hatglobal2and3}(a) shows the situation for $\beta_2=0.7$ with the points  $\mathbf{\widehat{B}_{2}^{\textit{h}}}$ and $\mathbf{\widehat{B}_{2}^{\ell}}$ that are created by $\mathbf{\widehat{CB}_{2}}$. At both $\mathbf{\widehat{B}_{2}^{\textit{h}}}$ and $\mathbf{\widehat{B}_{2}^{\ell}}$ the single loop $R^{*}$-symmetric periodic orbit has double Floquet multipliers $-1$, so this is a type of period-doubling bifurcation in this reversible context. Also shown here is the symmetry-breaking bifurcation $\mathbf{SB^{\textit{h}}}$ from Fig.~\ref{fig:Strans} that was created by $\mathbf{CSB}$ at $\beta_2=0.60052$. Note that the single loop $R^{*}$-symmetric periodic orbit at $\mathbf{SB^{\textit{h}}}$ has double Floquet multipliers $+1$, which constitutes the symmetry-breaking period-$k$ multiplying bifurcation with $k=1$. The branches of $R_1$-symmetric periodic orbits in Fig.~\ref{fig:hatglobal2and3}(a) that emerge from $\mathbf{SB^{\textit{h}}}$ now connect to the point $\mathbf{\widehat{B}_{2}^{\ell}}$ together with two further branches of $R_1$-symmetric periodic orbits. Together, they form the new surface ${}^{L}\mathcal{S}_{1}$ of $R_1$-symmetric periodic orbits with bounded period. The second new surface ${}^{M}\mathcal{S}_{1}$ consists of $R_1$-symmetric periodic orbits that emerge from $\mathbf{\widehat{B}_{2}^{\textit{h}}}$ and accumulate in a sprialling fashion on the $R_1$-symmetric primary homoclinic orbits. Note that ${}^{L}\mathcal{S}_{1}$} and ${}^{M}\mathcal{S}_{1}$ are invariant as surfaces under $R_1$ and $R_2$, but consists of $R_1$-symmetric periodic orbits; they are indeed the ``successors'' past $\mathbf{\widehat{CB}_{2}}$ of the surface  ${}^{B}\mathcal{S}_{1}$ from Fig.~\ref{fig:Sboth_trans} with the same property.

Figure~\ref{fig:hatglobal2and3}(b) is for $\beta_2=0.76$ and features also the points  $\mathbf{\widehat{B}_{2}^{\textit{h}}}$ and $\mathbf{\widehat{B}_{3}^{\ell}}$ that are created by $\mathbf{\widehat{CB}_{2}}$. In the process the new nested surfaces ${}^{N}\mathcal{S}_{1}$ and ${}^{O}\mathcal{S}_{1}$ are created. As we have seen before,  ${}^{N}\mathcal{S}_{1}$ connects $\mathbf{\widehat{B}_{3}^{\textit{h}}}$ and $\mathbf{\widehat{B}_{3}^{\ell}}$ on the surface ${}^{D}\mathcal{S}_{*}$ of $R^{*}$-symmetric periodic orbits, while ${}^{O}\mathcal{S}_{1}$ emerges from $\mathbf{\widehat{B}_{3}^{\textit{h}}}$ and spirals onto the $R_1$-symmetric primary homoclinic orbits. Indeed, this process of creating new nested surfaces continues for increasing $k$.

\section{Induced period-$k$ multiplying bifurcations in the zero-energy level}
\label{sec:order}

Periodic orbits in the zero-energy level of system~\eqref{eq:4Dode} are especially important: only they may form EtoP connections with $\mathbf{0}$, which, in turn, organise families of homoclinic orbits to $\mathbf{0}$, that is, solitons of the GNLSE \cite{PhysRevA.103.063514, PARKER2021132890}. This is why we now consider their induced bifurcations as $\beta_2$ is changed. 

As Fig.~\ref{fig:1parbif} already illustrated, $\Gamma^{*}$ becomes elliptic at the induced symmetry-breaking bifurcation $\mathbf{SB}$ at $\beta_2\approx0.60064$, where there is a double eigenvalue $+1$; see also Fig.~\ref{fig:Strans}(c2). Hence, for $\beta_2 > 0.60064$ the two non-trivial Floquet multipliers are a complex-conjugate pair of the form $\xi_{\pm}=e^{\pm 2\pi\alpha}$, where $\alpha$ is the rotation number of $\Gamma^{*}$. For increasing $\beta_2$ these Floquet multipliers move along the unit circle, namely $\xi_{+}$ counterclockwise and $\xi_{-}$ clockwise. In fact, both of them move along the entire unit circle and come back to form again a double eigenvalue $+1$ at $\beta_2\approx 0.8164$, which is the Hamiltonian-Hopf bifurcation $\mathbf{HH}$ (where the primary homoclinic orbits disappear). Hence, the rotation number $\alpha = \alpha(\beta_2)$ changes monotonically from $0$ at $\beta_2 = 0.60064$ to $1$ at $\beta_2 = 0.8164$. Whenever $\alpha(\beta_2)=\frac{p}{k}$ the Floquet multipliers $\xi_{\pm}$ are $k^{\text{th}}$ roots of unity and $\Gamma^{*}$ undergoes a period-$k$ multiplying bifurcation \cite{mackay1982renormalisation,Kent_1992}, which are points of $p\!:\!k$ resonance in the setting of reversible vector fields (or, equivalently, volume-preserving maps).

As we showed in Secs.~\ref{sec:periodmult} and ~\ref{sec:nonsymm}, the period-$k$ multiplying bifurcations points $\mathbf{{B}_{k}^{\ell}}$ and the symmetry-breaking period-$k$ multiplying bifurcations points $\mathbf{{\widehat{B}}_{k}^{\ell}}$ pass through the zero-energy level shortly after they are created at the points $\mathbf{{CB}_{k}}$ and $\mathbf{{\widehat{CB}}_{k}}$, respectively. In other words, there are two sequences of associated $\beta_2$-values when this happens and one finds the respective induced bifurcation in the zero-energy surface, which we refer to as $\mathbf{{B}_{k}}$ and $\mathbf{{\widehat{B}}_{k}}$, respectively. At at each such parameter point the single-loop $R^{*}$-symmetric periodic orbit in the zero-energy level undergoes the corresponding bifurcation, where $\beta_2$ is the bifurcation parameter. Note that this single-loop periodic orbit is the continuation in $\beta_2$ of the basic $R^{*}$-symmetric periodic orbit $\Gamma^{*}$ for $\beta_2 = 0.4$ from Sec.~\ref{sec:surfS*}; for notational simplicity, we refer to its continuation in the zero-energy level also as $\Gamma^{*}$ in this section.

\begin{figure}[t!]
   \centering
   \includegraphics[width=12cm]{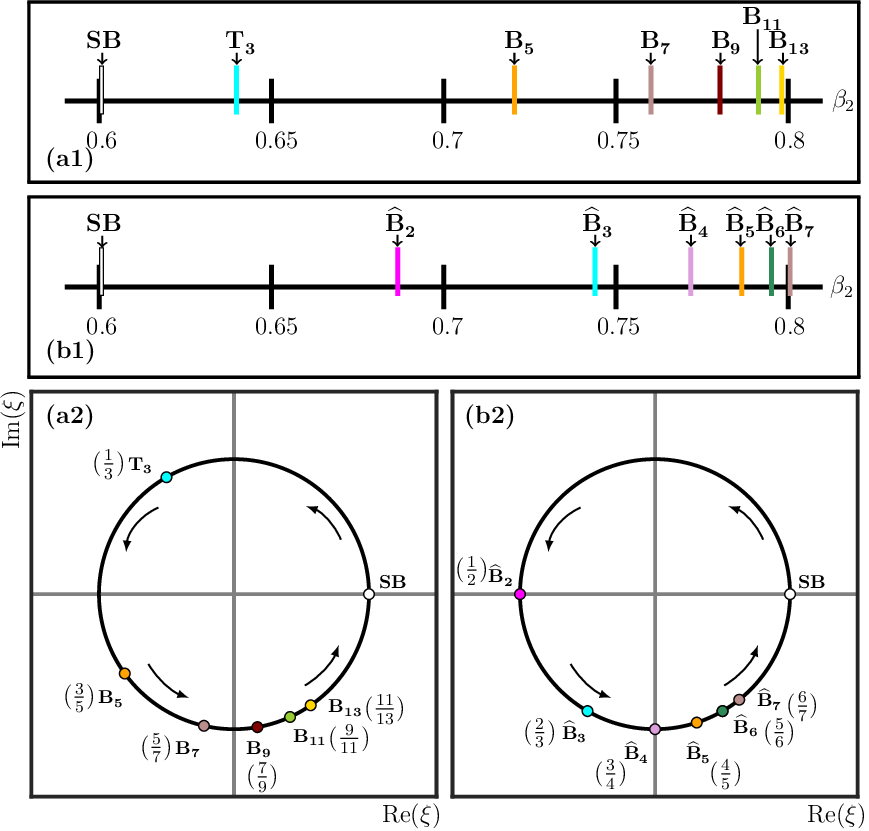}
       \caption{\label{fig:ordering} The ordering of $\beta_2$-values of the period-$k$ multiplying bifurcations of $\Gamma_*$ in the zero-energy level, and of the corresponding non-trivial Floquet multiplier $\xi_{+}$. The respective $\beta_2$-values and $\xi_{+}$ on the unit circle are shown in matching colours for $\mathbf{{B}_{k}}$ in panels~(a1) and~(b1) and for $\mathbf{{\widehat{B}}_{k}}$ in panels~(a2) and~(b2).} 
\end{figure} 

\begin{table}[h!]
\centering
\begin{tabular}{ |c|c|c|c|c|c|} 
\hline
\vtop{\hbox{\strut \vspace*{-2mm}}\hbox{\strut \hspace{0.2cm} period-$k$ }}& \vtop{\hbox{\strut}\vspace*{-2mm}\hbox{\strut$\beta_2$-value}} & \vtop{\hbox{\strut}\vspace*{-2mm}\hbox{\strut$\frac{p}{k}$}} & \vtop{\hbox{\strut symmetry-breaking}\hbox{\strut \hspace{7mm} period-$k$ }} & \vtop{\hbox{\strut}\vspace*{-2mm}\hbox{\strut$\beta_2$-value}} & \vtop{\hbox{\strut}\vspace*{-2mm}\hbox{\strut$\frac{p}{k}$}} \\ [2pt]
\hline 
& & & & & \\[-8pt]
\quad $\mathbf{SB}$ ($= \mathbf{{B}_{1}}$) & $0.60064$  & $\frac{-1}{1}$& \quad $\mathbf{SB}$ ($= \mathbf{\widehat{B}_{1}}$)& $0.60064$ & $\frac{0}{1}$\\ [2pt] 
\quad $\mathbf{{T}_{3}}$ ($= \mathbf{{B}_{3}}$)& $0.639829$  & $\frac{1}{3}$& $\mathbf{\widehat{B}_{2}}$ & 0.68667 & $\frac{1}{2}$\\ [2pt] 
$\mathbf{{B}_{5}}$ & $0.720559$ & $\frac{3}{5}$& $\mathbf{\widehat{B}_{3}}$ & 0.743907 & $\frac{2}{3}$\\ [2pt]
$\mathbf{{B}_{7}}$ & $0.760154$ & $\frac{5}{7}$& $\mathbf{\widehat{B}_{4}}$ & 0.771728 & $\frac{3}{4}$\\ [2pt]
$\mathbf{{B}_{9}}$ & $0.780195$ & $\frac{7}{9}$& $\mathbf{\widehat{B}_{5}}$ & 0.786516 & $\frac{4}{5}$\\ [2pt]
$\mathbf{{B}_{11}}$ & $0.791361$ & $\frac{9}{11}$& $\mathbf{\widehat{B}_{6}}$ & 0.795141 & $\frac{5}{6}$\\ [2pt]
$\mathbf{{B}_{13}}$ & $0.798140$ &  $\frac{11}{13}$& $\mathbf{\widehat{B}_{7}}$ & 0.800559 & $\frac{6}{7}$\\ [2pt]
 .&.  &  .&. &. &.  \\[2pt]
  .&.  &  .&. &. &.  \\[2pt]
 .&.  &  .&.  &. &.  \\[2pt]

$\mathbf{{B}_{k}}$ & for odd $k\geq3$  & $\frac{k-2}{k}$ & $\mathbf{\widehat{B}_{k}}$ & for all $k\geq2$ & $\frac{k-1}{k}$\\ [2pt]

  .&.  &  .&. &. &.  \\[2pt]
 .&.  &  .&.  &. &.  \\[2pt]
 
\vtop{\hbox{\strut \hspace{1.1cm} $\mathbf{HH}$}\hbox{\strut \hspace{0.25cm}  for odd $ k \to \infty$ }}& $0.8164$  & $1$ & \vtop{\hbox{\strut \hspace{0.7cm} $\mathbf{HH}$}\hbox{\strut \hspace{0.25cm}  for $ k \to \infty$ }} & $0.8164$  & $1$\\[15pt]
\hline
\end{tabular} \vspace*{3mm}
\caption{\label{table:ordering} The $\beta_2$-values and corresponding rational rotation numbers of the computed period-$k$ multiplying bifurcations $\mathbf{{B}_{k}}$ and $\mathbf{{\widehat{B}}_{k}}$ of $\Gamma_*$ in the zero-energy level.}
\end{table}

Figure~\ref{fig:ordering} illustrates that the two sequences $\mathbf{{B}_{k}}$ and $\mathbf{{\widehat{B}}_{k}}$ of different types of period-$k$ multiplying bifurcations of $\Gamma^{*}$ occur at specific subsequences of roots of unity. Panels~(a1) and~(b1) show the $\beta_2$-values where $\mathbf{{B}_{k}}$ and $\mathbf{\widehat{B}_{k}}$ take place, respectively. Panels~(a2) and~(b2) show the corresponding positions of the positive nontrivial Floquet multiplier $\xi_{+}$ of $\Gamma_{*}$ on the unit circle in the complex plane; here
the rational rotation numbers are given and the arrows indicate that $\xi_{+}$ moves counterclockwise for increasing $\beta_2\in [0.60064, 0.8164]$. The $\beta_2$-values and rational rotation numbers at the computed points $\mathbf{{B}_{k}^{\ell}}$ and $\mathbf{\widehat{B}_{k}^{\ell}}$ were identified during the continuation of $\Gamma^{*}$. They are listed in Table~\ref{table:ordering}, where we make the identifications $\mathbf{SB} = \mathbf{{B}_{1}}$, $\mathbf{{B}_{3}} = \mathbf{{T}_{3}}$, and $\mathbf{SB} = \mathbf{\widehat{B}_{1}}$; note also that the rotation number is taken modulo $1$, meaning that $\frac{-1}{1} = \frac{0}{1} = \frac{1}{1}$. Hence, the $\beta_2$-values of the sequences both start at $\mathbf{SB}$ with $0.60064$, and they are interleaved on the $\beta_2$-line; see Fig.~\ref{fig:ordering}(a1) and~(b1). The data we computed allows us to draw the following conclusions, which we formulate as: 

\begin{conjecture}
\label{conj:ordering}
Ordering of the bifurcations $\mathbf{{B}_{k}}$ and $\mathbf{{\widehat{B}}_{k}}$ of the single-loop $R^{*}$-symmetric periodic orbit $\Gamma^{*}$ in the zero-energy level.
\begin{enumerate} 
\item The period-$k$ multiplying bifurcations  $\mathbf{{B}_{k}}$ occur at the values of $\beta_2$ where $\xi_{\pm}=e^{\pm 2\pi\frac{k-2}{k}}$ for all odd $k \geq 1$. 
\item The symmetry-breaking period-$k$ multiplying bifurcations $\mathbf{{\widehat{B}}_{k}}$ occur at the values of $\beta_2$ where $\xi_{\pm}=e^{\pm 2\pi\frac{k-1}{k}}$ for all odd $k \geq 1$. 
\item 
The corresponding bifurcating $R^{*}$-symmetric, $R_{1}$-symmetric and $R_{2}$-symmetric $k$-loop orbits exist below and up the respective $\beta_2$-value; except for the transcritical-type case $\mathbf{{B}_{3}} = \mathbf{{T}_{3}}$, where the triple-loop orbits exist up $\mathbf{SN}$ at $\beta_2 \approx 0.649089$.
\item Since $\frac{k-2}{k}$ and $\frac{k-1}{k}$ tend to $1$ as $k \to \infty$, both sequences $\mathbf{{B}_{k}}$ and $\mathbf{\widehat{B}_{k}}$ approach the point $\mathbf{HH}$ and, hence, the associated value $\beta_2 \approx 0.8164$.
\end{enumerate}
\end{conjecture} 

The observations presented in Conjecture~\ref{conj:ordering} are supported by the numerical results presented in this paper as summarised in Table~\ref{table:ordering}. The upshot is that $k$-loop periodic orbits of different symmetry type disappear from the zero-energy level one-by-one at the interleaving points of the two sequences $\mathbf{{B}_{k}}$ and $\mathbf{\widehat{B}_{k}}$ as $\beta_2 \approx 0.8164$ is increased towards the point $\mathbf{HH}$. Indeed, multi-loop periodic orbits do no longer exist in the zero-energy level at $\mathbf{HH}$, where the primary homoclinic orbits disappear.

\section{Conclusions}
\label{sec:conclusions}

We presented a detailed study of the underlying periodic orbit structure of the four-dimensional Hamiltonian system \eqref{eq:4Dode} with two reversible symmetries, which arises from a traveling-wave ansatz for solitons in the GNLSE \eqref{eq:gnlse} with quadratic and quartic dispersion. We focused on three basic basic surfaces $\mathcal{S}_{*}$, $\mathcal{S}_{1}^{+}$ and $\mathcal{S}_{1}^{-}$ and $\mathcal{S}_{*}$, which are directly associated with the basic $R_1$-symmetric homoclinic orbit and its $R_2$-symmetric counterpart. These three surfaces are initially distinct from one another when the quadratic dispersion parameter $\beta_2$ is not too big. However, they interact with each other in complicated ways via different bifurcations for increasing $\beta_2$, and we illustrated how the surfaces of periodic orbits and their intersection sets with a three-dimensional section change as a result. Specifically, we identified symmetry-breaking bifurcations, period-$k$ multiplying bifurcations, and saddle-node bifurcations. Each of these bifurcations has a degenerate case at a specific value of $\beta_2$, which always takes place at a particular energy level $H>0$. The degenerate bifurcations change the overall geometry of the surfaces of periodic orbits by splitting certain surfaces and creating new ones, which join at a created pair of period-$k$ multiplying bifurcations.  We illustrated in $(u_1, u_2, H)$-space how an increasing number of nested surfaces of periodic orbits of different symmetry properties as $\beta_2$ is increased. In the zero-energy level, this scenario is reflected by associated sequences of induced $k$-multiplying bifurcations at which the (continuation of the) basic single-loop periodic orbit $\Gamma_{*}$ has Floquet multipliers from specific subsequences of roots of unity. We formulated these observations as Conjecture~\ref{conj:ordering}. 

The comprehensive description of what is happening to the three basic surfaces of periodic orbits when $\beta_2$ is increased towards the Hamiltonian-Hopf bifurcation \textbf{HH} may serve as a blueprint for understanding the overall structure of all periodic orbits of system \eqref{eq:4Dode}, which is actually even more complicated. To begin with, the Floquet multipliers of $\Gamma^{*}$ cross many more roots of unity, so there will be infinitely many other sequences where multi-loop periodic orbits disappear from the zero-energy level. Moreover, the nontrivial complex conjugate Floquet multipliers of elliptic multi-loop periodic orbits also go through sequences of resonances, where further branches of higher-period periodic orbits bifurcate. The picture that emerges as an extension of Conjecture~\ref{conj:ordering} is one of infinitely many sequences of certain types of periodic orbits disappearing from the zero-energy level as $\beta_2$ is increased towards $\mathbf{HH}$ --- forming a dense set of $\beta_2$-values where the corresponding Floquet multipliers are roots of unity. Moreover, in light of the rescaling result in \cite{PhysRevA.103.063514}, each of these $k$-multiplying bifurcations will form a parabola in the $(\beta_2,\mu)$-plane shown in Fig.~\ref{fig:connections}(a).

It is a natural conjecture that the underlying mechanism for this scenario in the zero-energy level is again the one we found here: the emergence of special degenerate periodic orbits at certain values of $\beta_2$ that create pairs of $k$-multiplying periodic orbits off the zero-energy level, one of which then moves through the zero-energy level to induce the corresponding bifurcation there. Our results strongly suggest that such degenerate periodic orbits, of which the sequences $\mathbf{{CB}_{k}}$ and $\mathbf{\widehat{CB}_{k}}$ are examples, exist at a set of $\beta_2$-values that is dense in the interval $[0.60064, 0.8164]$, that is, between $\mathbf{SB}$ and $\mathbf{HH}$. Moreover, we conjecture that this set is actually dense in the larger interval $[-0.8164, 0.8164]$ --- hence, all the way down to the Belyakov-Devaney bifurcation \textbf{BD}. Again, these degenerate bifurcations will be found along parabolas in the $(\beta_2,\mu)$-plane.

\section{Acknowledgements}

R. Bandara thanks the Dodd-Walls Centre for Photonic and Quantum Technologies for support of this research by a PhD scholarship held at the University of Auckland.
A. Giraldo was supported at the Korea Institute for Advanced Study by KIAS Individual Grant No. CG086101.

\providecommand{\noopsort}[1]{}\providecommand{\singleletter}[1]{#1}%

\end{document}